\documentclass[12pt, twoside, epsf]{article}

\usepackage{color,graphicx,times}
\usepackage{amsmath, amssymb, graphics}

\newcommand{\mathsym}[1]{{}}
\newcommand{\unicode}[1]{{}}

\let \ttorg \tt \def \tt{\ttorg \obeyspaces}

\begin{document}

 \date{}

\title{\bf Idempotents in the Temperley-Lieb Monoid and Other Categories}

\author{Louis H. Kauffman \\
  Department of Mathematics, Statistics and Computer Science \\
  University of Illinois at Chicago \\
  851 South Morgan Street\\
  Chicago, IL, 60607-7045}

\maketitle
  
\thispagestyle{empty}

\begin{abstract}
 This paper examines idempotents in algebras and categories that arise from factorizations of the identity morphism. In the diagrammatic and combinatorial contexts considered here, these factorizations correspond to 
 generalizations of meanders, where a meander is understood to be a curve in the plane that wanders transversely back and  forth across a given straight line. By formulating the algebra of the Temperley-Lieb Monoid in terms of 
 planar curve combinatorics, one can understand idempotents in the Temperley-Lieb Monoid in terms of meanders. Corresponding results are shown for the Brauer Monoid and for the Tangle Monoid and Tangle Category.\\

\end{abstract}

\bigbreak



\section{Introduction}
\begin{figure}[htb]
     \begin{center}
     \begin{tabular}{c}
     \includegraphics[width=8cm]{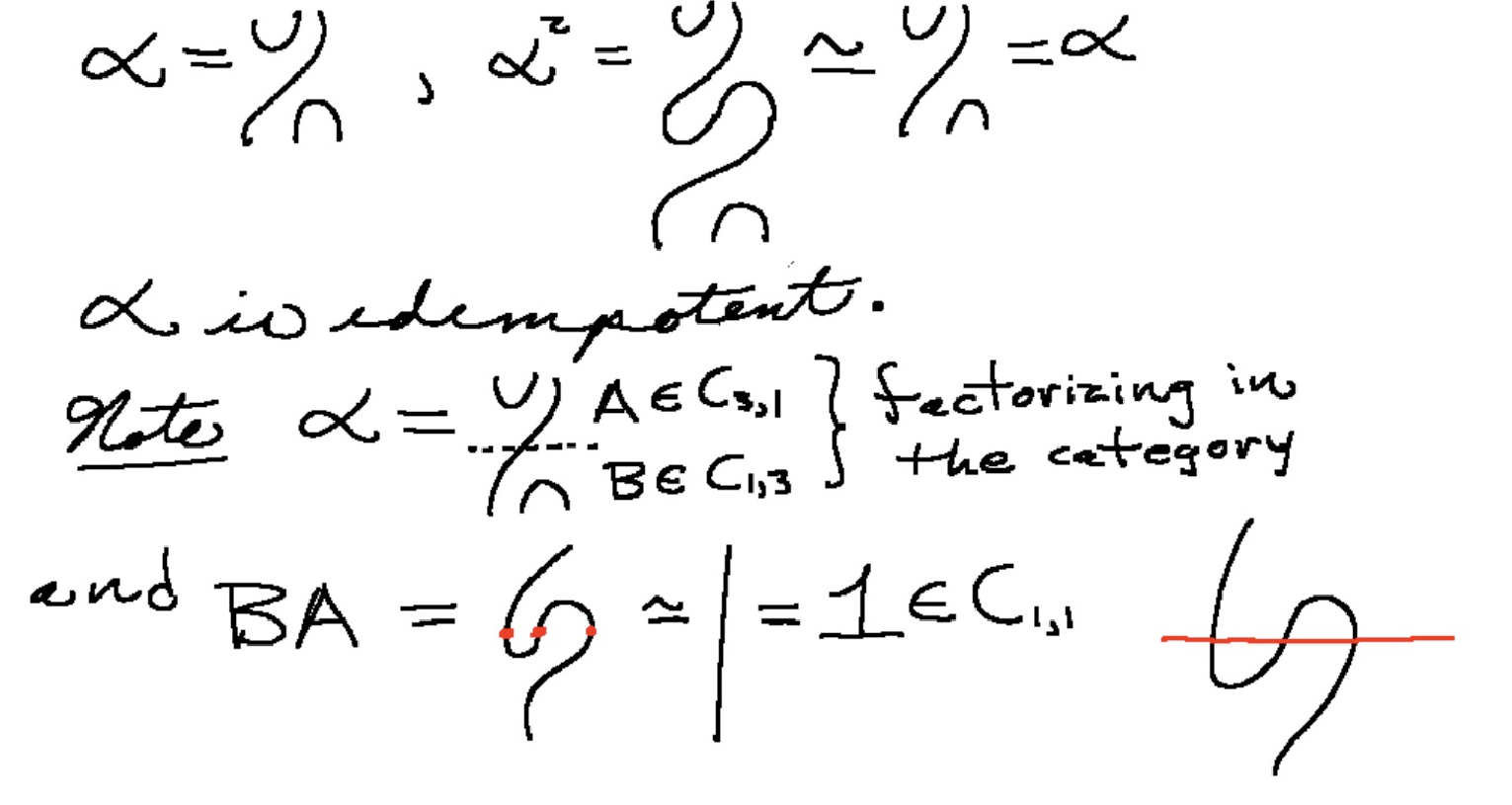}
     \end{tabular}
     \caption{\bf An idempotent constructed by factorization of the identity}
     \label{examples}
\end{center}
\end{figure}

This paper is about the construction and properties of idempotents, elements $P$ in an algebra such that $P^2=P$, in many algebraic situations.
To introduce these ideas, we begin with the example shown in Figure~\ref{examples}. In the figure we see an element $\alpha$ of the Temperley-Lieb
Connection Monoid (definitions will occur in later sections of the paper) and a diagrammatic illustration showing that $\alpha^2 = \alpha.$ 
What is it in the structure of $\alpha$ that makes it an idempotent? The second part of the figure illustrates that in the Connection Monoid Category, $\alpha$ factorizes as
$\alpha = AB$ where $A$ is a $[3]$ to $[1]$ morphism while $B$ is a $[1]$ to $[3]$ morphism. The composition $BA$ is well-defined and is seen to be a factorization of the identity 
from $[1]$ to $[1]$. Thus $BA = 1$ and we see that $$\alpha^2 = ABAB = A(BA)B = A1B = AB = \alpha.$$ This method can be generalized to produce infinitely many idempotents in the Temperley-Lieb
Monoid, and indeed all idempotents if we allow multiple strands in the factorization.\\

If $Cat$ is any category, and $BA = 1$ is a factorization of the identity in that category, then $P = AB$ will be an idempotent morphism in the category. The proof is the same as above.
$$P^2 = ABAB = A(BA)B = A1B = AB = P.$$ Note that composition in a category is associative.\\

In this paper we study idempotents in the Temperley-Lieb Monoid and the structure of powers of elements in the monoid. Then we will apply these ideas to other categories showing how idempotents arise in the Brauer Monoid and in the Tangle Category and Tangle Monoid.  It is of interest also to see how this idea of factorizing the identity in a category can be related to other structures such as the replication of DNA. See \cite{BL1,BL2,AP,KLogic,MLogic}. The idea behind self-replication in this mode is very simple. Given that $BA=1$, we can consider a process that goes in the opposite order to reducing an idempotent:$$AB \longrightarrow A1B \longrightarrow ABAB \longrightarrow AB AB.$$ One can interpret the separation of $AB$ into $A1B$ as analogous to the separation of the Watson and Crick strands by enzyme action in the cell. Then $1$ stands for the environment and $1 \longrightarrow B A$ corresponds to the way that the environment supplies complementary base pairs to the separated strands. This idea has been pursued in the references just cited and will be further examined in a sequel to the present work.\\

\begin{figure}[htb]
     \begin{center}
     \begin{tabular}{c}
     \includegraphics[width=9cm]{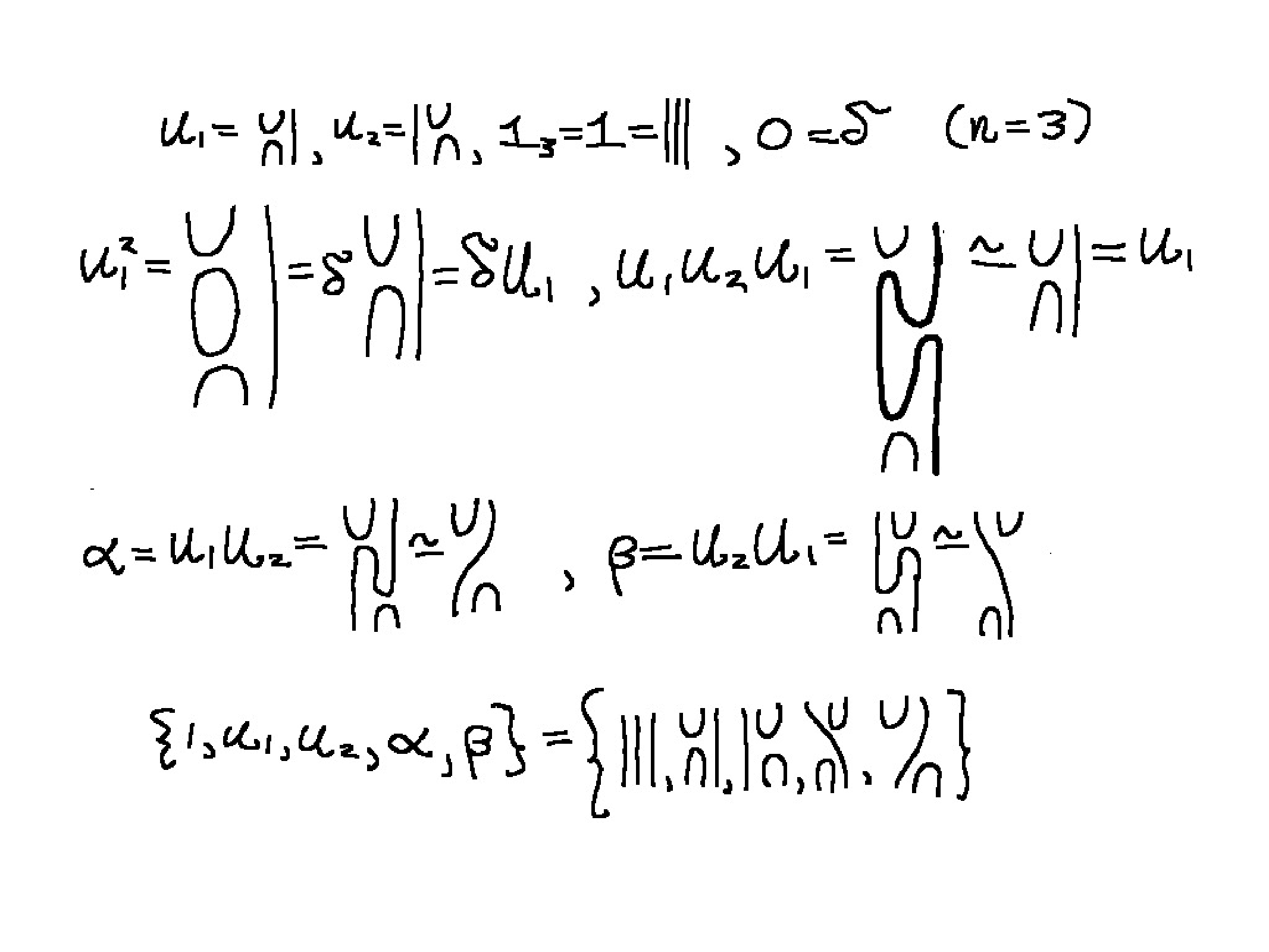}
     \end{tabular}
     \caption{\bf Diagrammatic Temperley-Lieb relations}
     \label{TLDiagrams}
\end{center}
\end{figure}

\begin{figure}[htb]
     \begin{center}
     \begin{tabular}{c}
     \includegraphics[width=7cm]{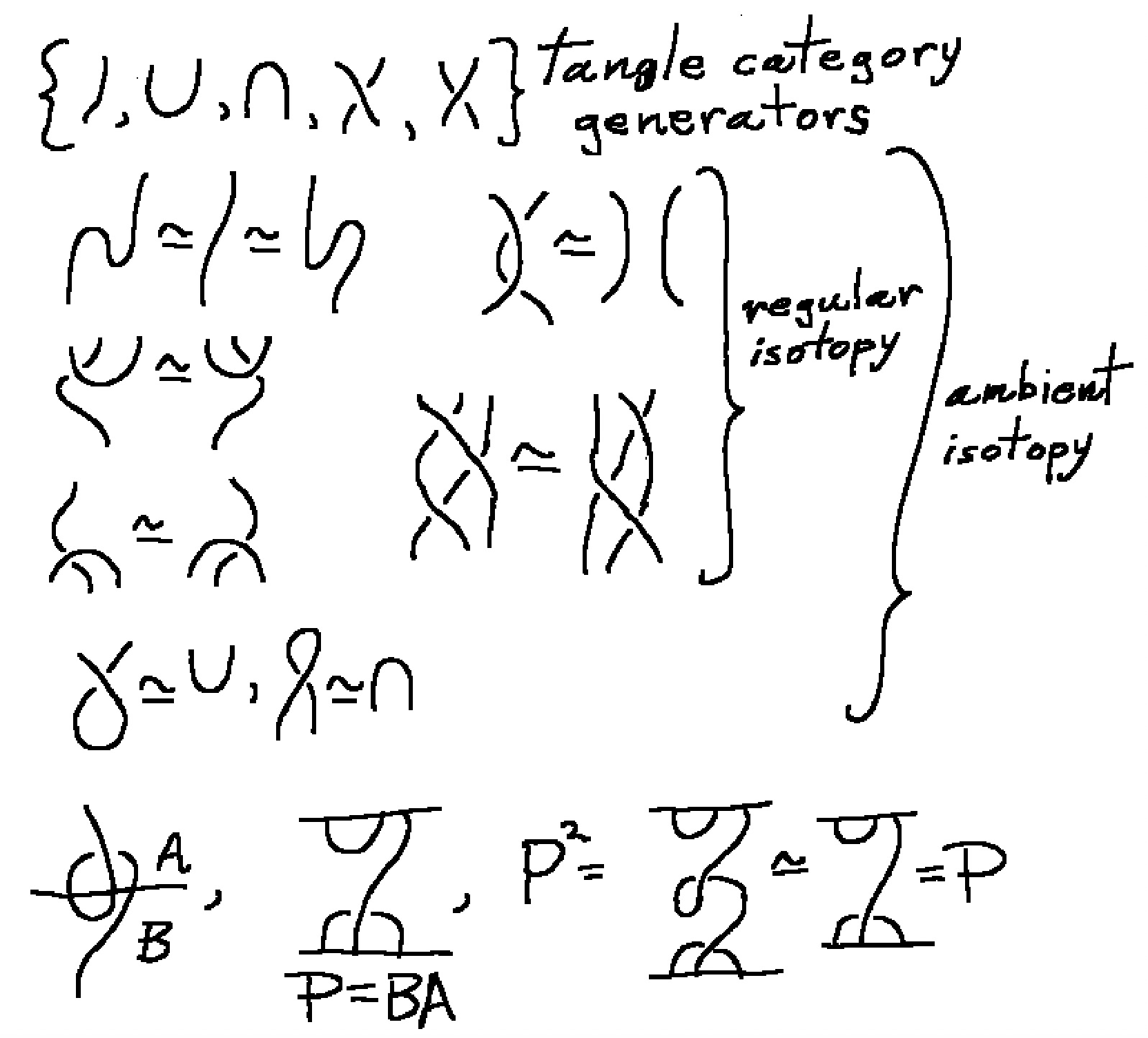}
     \end{tabular}
     \caption{\bf The Tangle Category}
     \label{tanglecat}
\end{center}
\end{figure}

\begin{figure}[htb]
     \begin{center}
     \begin{tabular}{c}
     \includegraphics[width=8cm]{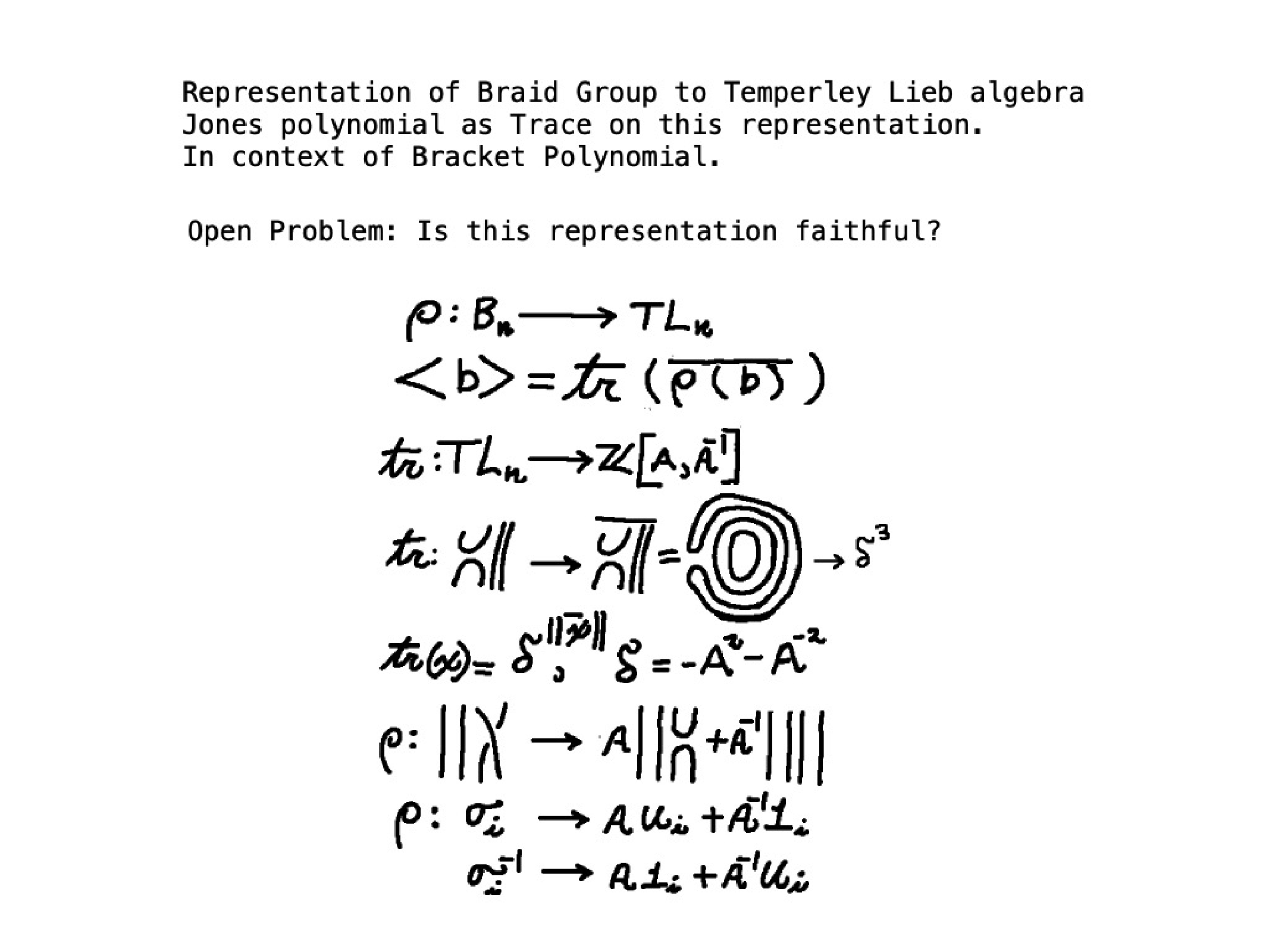}
     \end{tabular}
     \caption{\bf Representation Question}
     \label{repques}
\end{center}
\end{figure}

\begin{figure}[htb]
     \begin{center}
     \begin{tabular}{c}
     \includegraphics[width=8cm]{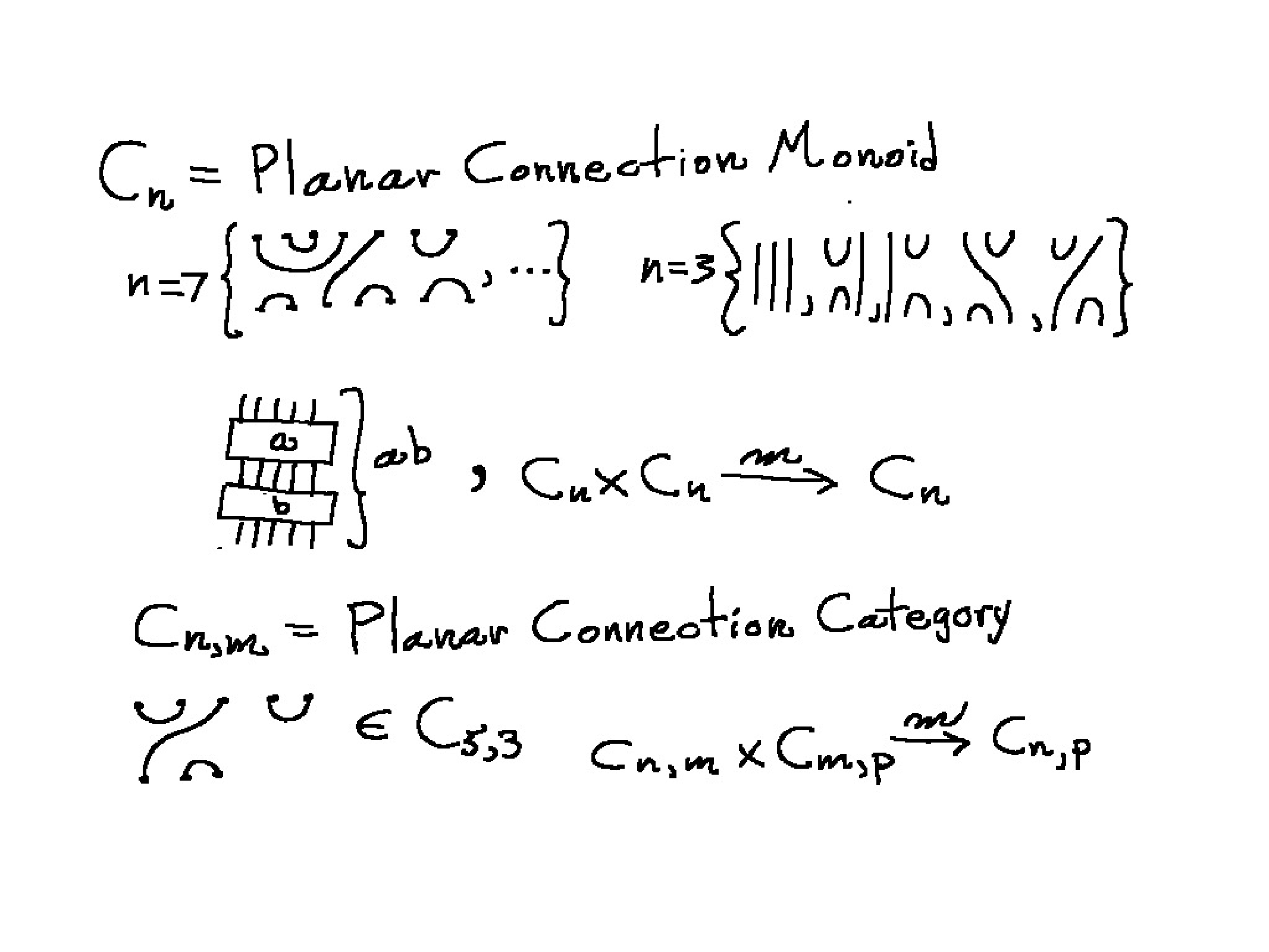}
     \end{tabular}
     \caption{\bf The Planar Connection Monoid}
     \label{conncat}
\end{center}
\end{figure}

\begin{figure}[htb]
     \begin{center}
     \begin{tabular}{c}
     \includegraphics[width=7cm]{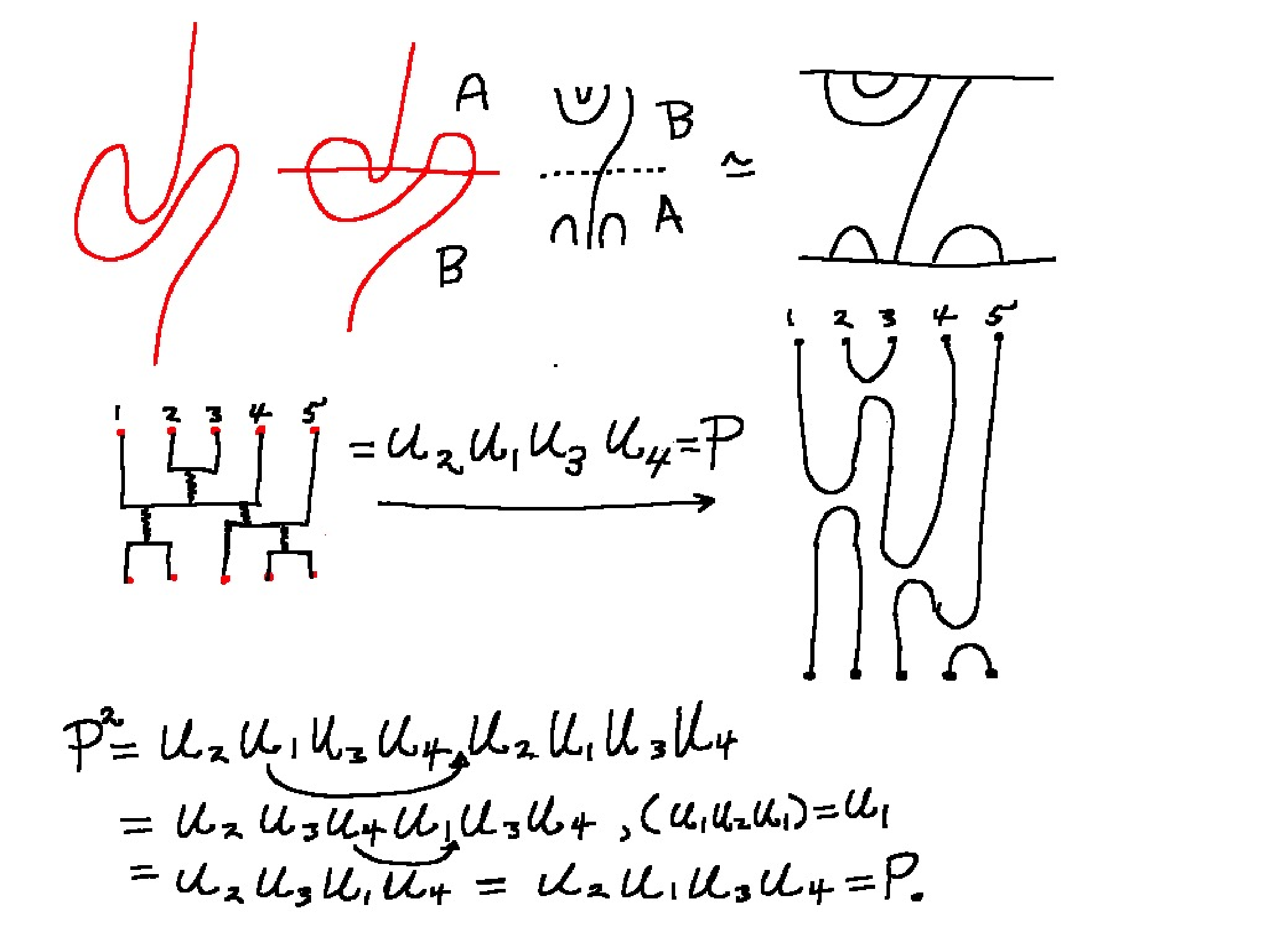}
     \end{tabular}
     \caption{\bf Rectilinear Method converting the Connection Monoid to the Temperley-Lieb Monoid}
     \label{rect1}
\end{center}
\end{figure}

\begin{figure}[htb]
     \begin{center}
     \begin{tabular}{c}
     \includegraphics[width=7cm]{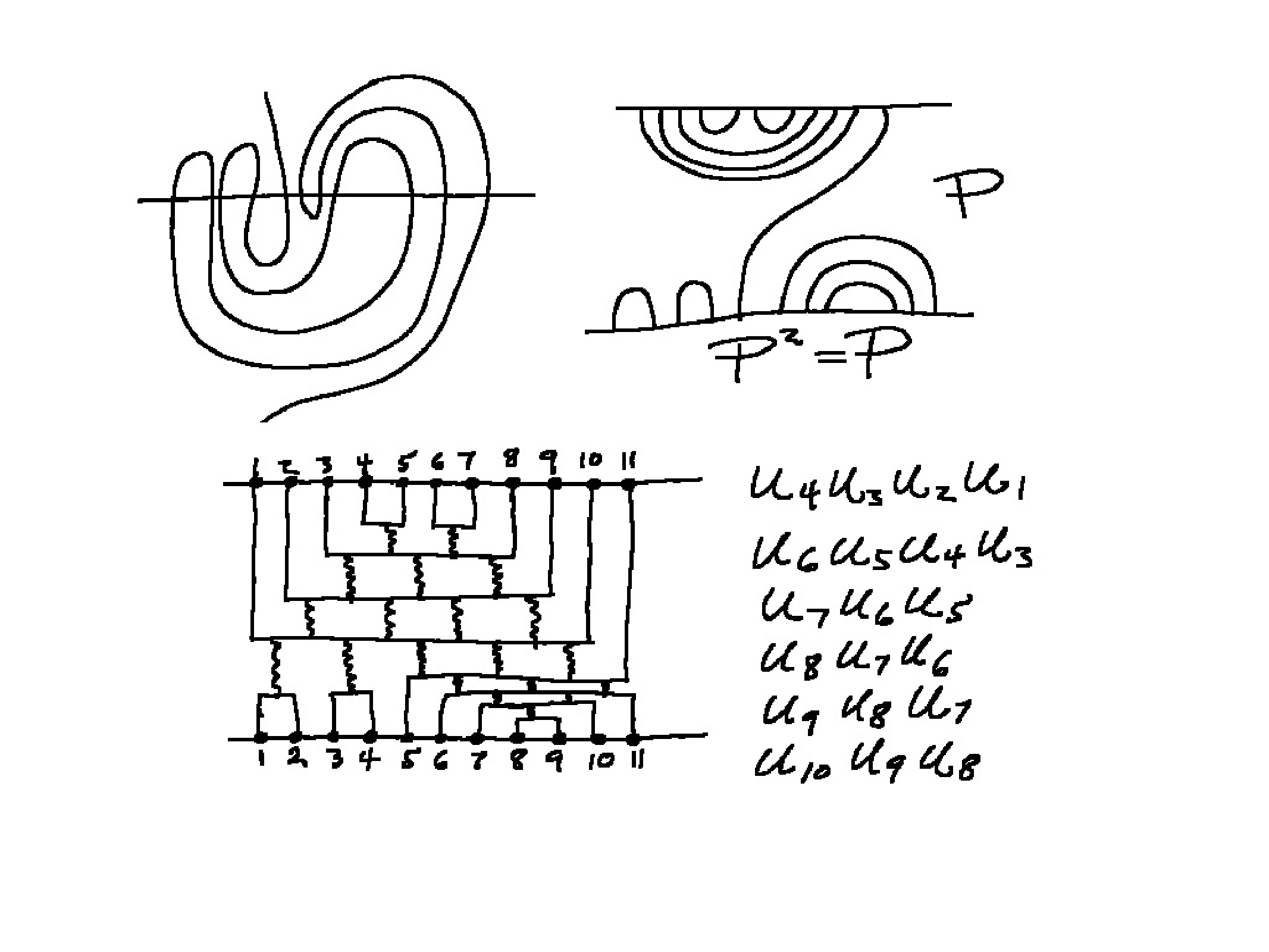}
     \end{tabular}
     \caption{\bf Example of Rectilinear Conversion}
     \label{rect2}
\end{center}
\end{figure}

\section{The Brauer and Temperley-Lieb Categories}
In this section we first describe the Temperley-Lieb Monoid ($TLM_n$) (on $n$ strands) \cite{KD,Bracket,Stat,KP,KLins} in terms of generators and relations and then reformulate it in terms of the Connection Monoid ($CM_n$) with its diagrammatic interpretation in terms of arcs connecting two rows of $n$ points in the plane. We generalize the Connection Monoid to a category by allowing the two rows to have different numbers of points. It is often useful to factor elements of the Connection Monoid into compositions in the Connection Category. Thus we can speak of the Temperley-Lieb Category as another name for the Connection Category. The element $P$ in Figure~\ref{examples} is factored into a product of elements in the Connection Category.\\

A word about categories is in order. In Figure~\ref{tanglecat} we illustrate the generators and basic relations for morphisms in the {\it Tangle Category} (some relations are left out of the diagram such as the commutation of cups and caps). This category stands behind many of the constructions in this paper. It is a tensor category, meaning that the basic morphisms can be both composed by attaching lines vertically, and juxtaposed by placing them horizontally. Morphisms are arranged so that they have cups and caps and crossings according to a height function on the corresponding diagram and so are sometimes called Morse diagrams. The relations in the category consist in braiding identities (a subset of the Reidemeister moves), cancellations of cups and caps, and swing moves that shift a crossing over a cup or a cap to the opposite crossing type. If we remove crossings from the Tangle Category we obtain the Connection Category. If we flatten the crossings in the Tangle Category, we obtain the Brauer Category discussed below. The algebraic structures such as the Temperley-Lieb Monoid and the Brauer Monoid described below have algebraic relations that are related to these moves in the category. See \cite{KD} for further discussion of these relationships.\\ 

By using factorizations of the identity in the Connection Category we construct idempotents in the Temperley-Lieb Monoid. These factorizations are identified with {\it meanders} \cite{Meanders} in the diagrammatic representation of the Connection Category. We show how to articulate these idempotents in the Temperely-Lieb Monoid in terms of the generators for that monoid. It is an interesting problem to characterize these idempotents at the purely algebraic level where the monoid is described in terms of generators and relations. We then continue in the following section with a study of powers of elements in the Temperley-Lieb Monoid, using these techniques and with descriptions of numerous examples.\\

\noindent {\bf Definition.} The Temperley-Lieb Monoid ($TLM_n$) is an associative algebraic structure with one binary relation, defined relative to an element $\delta$ with no relations other than that $\delta$ commutes with all elements in the monoid. The monoid is generated by $U_1, U_2, \cdots U_{n-1}$ with the relations
\begin{enumerate}
\item $U_i^2 = \delta U_i$ for all $i = 1, 2, \cdots , n-1.$
\item $U_i U_j U_i = U_i$ whenever $|i-j|=1.$
\item $U_i U_j = U_j U_i$ whenever $|i-j|>1.$
\end{enumerate}

\begin{figure}[htb]
     \begin{center}
     \begin{tabular}{c}
     \includegraphics[width=7cm]{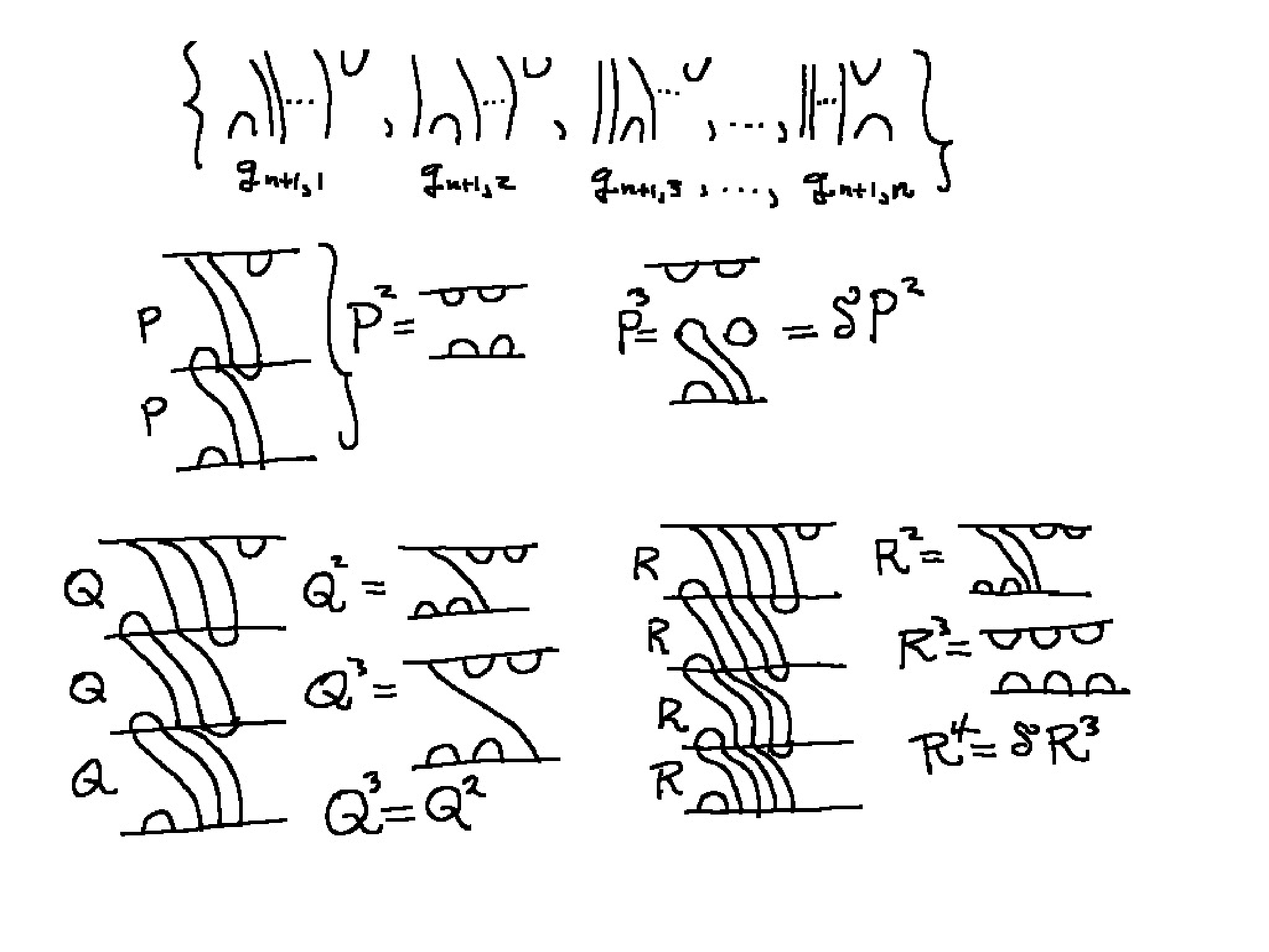}
     \end{tabular}
     \caption{\bf Powers Example 3}
     \label{example3}
\end{center}
\end{figure}

\noindent {\bf Remark.} The {\it Temperley-Lieb Algebra $TL_{n}$ } is defined by letting $U_1, U_2, \cdots U_{n-1}$ generate a module over the ring $R = Z[\delta]$ so that one can add elements of the monoid freely to one another with coefficients in $R.$
For many applications, idempotents in this algebra (such as the Jones-Wenzl  projectors \cite{KLins,Morrison}) constructed using both the additive and multiplicative structure in the algebra are of great interest. In a subsequent paper we will discuss the implications of our analysis of monoid idempotents and powers of monoid elements for the structure of these more general idempotents. A beginning for this consideration occurs in our Figure~\ref{example3} where we show how powers of elements behave that are used by Scott Morrison \cite{Morrison} to construct Jones-Wenzl projectors.\\

Before introducing the Connection Monoid, note the diagrams in Figure~\ref{TLDiagrams}. In these diagrams $U_i$ is represented by arcs connecting two rows of $n$ points so that the $k$-th point on the upper row is connected to the $k-th$ point on the lower row for $k \ne i$ and $k \ne i+1,$ and on the top row $i$ is connected to $i+1$ and on the bottom row the $i$-th point is connected to the $i+1$-th point. The diagram for $U_i$ can be described as having a paired minimum (from the top row) and maximum (from the bottom row. These minima and maxima are the result of connecting smooth arcs from top to top and from bottom to bottom. We then see, taking the value of a loop to be $\delta$, that all the relations in the Temperley-Lieb Monoid are satisfied by taking these diagrams up to topological deformation in the plane and by composing the diagrams A and B  by joining the bottom row of points on A to the top row of points on B to form a new planar diagram. Note in the figure how the topological deformation of the arcs results in the relations in the monoid and that when we multiply $U_i$ by itself, a loop forms in the middle, resulting in the first relation for the monoid. These diagrams then generate an image of the Temperley-Lieb Monoid. \\

Each diagram in Figure~\ref{TLDiagrams} consists in a choice of connections among the points of the top and bottom rows such that the entire set of connecting arcs can be embedded in the plane space between the top and bottom rows of points without any intersections between arcs or self-intersections. In the Temperley-Lieb diagrams we have also written them so that they can be directly interpreted (via the paired maxima and minima) as products of the generators $U_i.$ \\

\begin{figure}[htb]
     \begin{center}
     \begin{tabular}{c}
     \includegraphics[width=7cm]{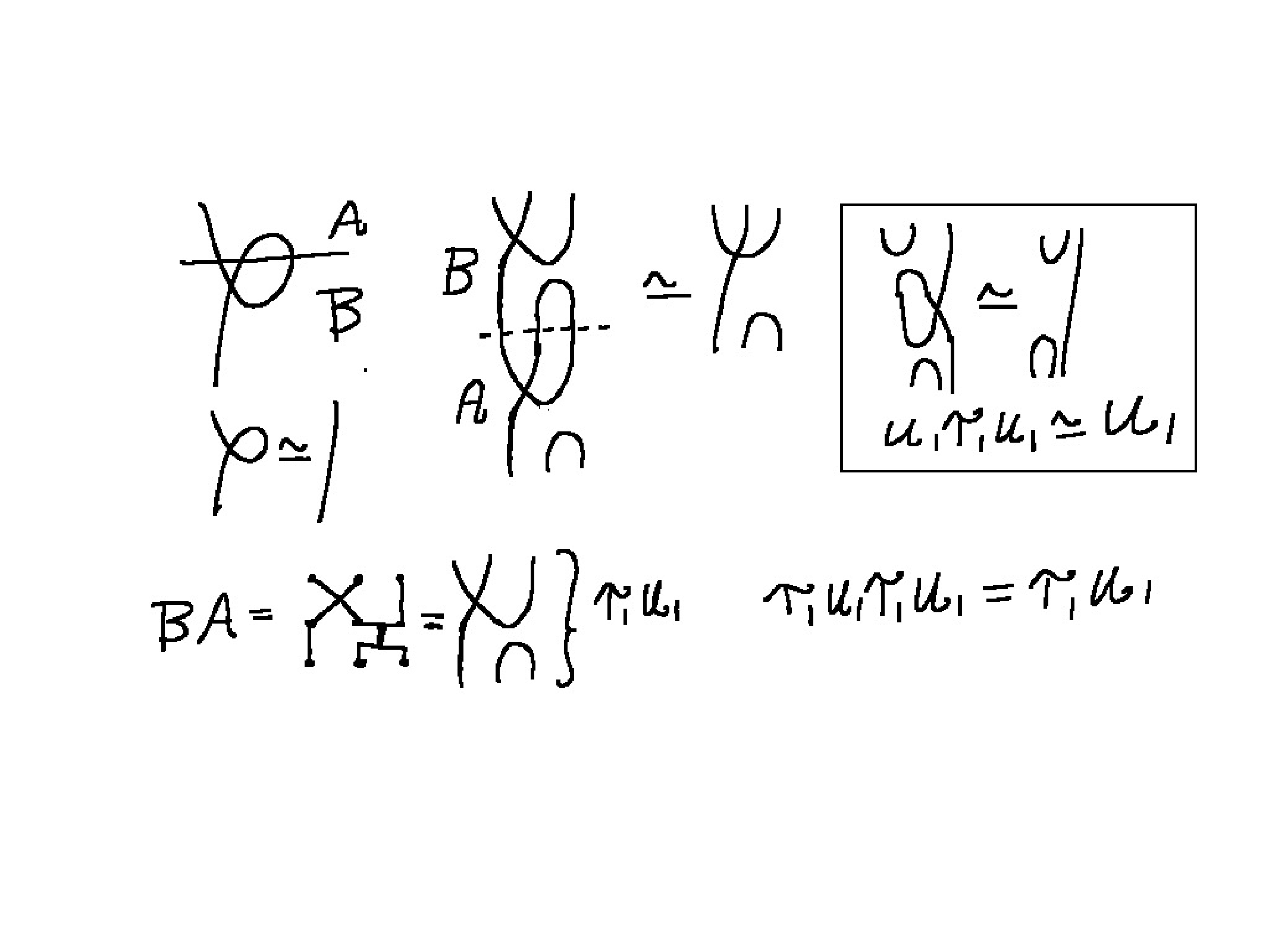}
     \end{tabular}
     \caption{\bf Constructing an Idempotent in the Brauer Monoid}
     \label{Brauer1}
\end{center}
\end{figure}

\noindent {\bf Definition.} The {\it Brauer Connection Monoid} $BCM_{n}$ consists in two collections of nodes $1,2,\cdots ,n$ (the {\it top row}) and $1',2',\cdots ,n'$ (the {\it bottom row}) and all pairings among these nodes so that all the nodes are partitioned into pairs of the form $(ij), (ij'), (i'j')$ indicating pairings from top row to top row, top row to bottom row and bottom row to bottom row. For example, if $n=4$ then $a = (12')(24)(31')(3'4')$ is an example of a pairing. Two pairings $a$ and $b$ are multiplied by ``attaching the bottom row of $a$ to the top row of $b$". We accomplish this so that $(ij') \times (jk') = (ik')$, that is the matching nodes are eliminated to create a new pairing. Elements of the Brauer Connection Monoid can be instantiated as diagrams by using arcs to indicate the connections in the space between the top and bottom rows. See Figure~\ref{Brauer1} for an illustration of elements in the diagrammatic Brauer Monoid. Note that in these diagrams, arcs can cross one another. The collection of arcs for a given element of the Brauer Connection Monoid are taken up to topological deformation in the space between the top and bottom rows. The figure illustrates equivalences and particularly equivalences after multiplications. See \cite{Brauer} for more information about the Brauer Monoid and its history. Note that the Brauer Connection Monoid is
defined independent of the diagrams that represent it and that if we take two diagrams to be equivalent whenever they represent the same pairings, then the diagrammatic is a faithful representation of the Monoid.\\

\noindent {\bf Definition.} The {\it Connection Monoid} $CM_{n}$ is the sub-monoid of the Brauer Connection Monoid consisting in all pairings that can be represented by collections of embedded arcs with no crossings between any pair of arcs and no self-crossings of any arc.
Thus each element in the diagrammatic Temperley-Lieb Monoid can be seen as an element of the Connection Monoid. Just as in the Brauer Connection Monoid, two elements of the Connection Monoid are equal exactly when they make identical sets of pairs. Connection Monoid diagrams are taken up to topological equivalence as in the Brauer Connection Monoid. Note that the Brauer Connection Monoid $BCM_{n}$ includes a copy of the symmetric group $S_{n}$ of permutations of $n$ nodes, while these permutations, involving crossing arcs, are missing from the Connection Monoid.\\

\noindent {\bf Theorem.} {\it  For every $n,$ the Temperley-Lieb Monoid $TLM_n$ and the Connection Monoid $CM_n$ are isomorphic.} \\

\noindent {\bf Proof.} Given any element of the Connection Monoid, we associate a specific product of generators in the Temperley-Lieb Monoid. View Figure~\ref{rect1} and Figure~\ref{rect2}. These figures illustrate how to draw a {\it rectilinear} diagram for any element $Q$ of $CM_{n}.$
Each connecting arc is a composite of horizontal and vertical lines. Any connection from the top row to the bottom row consists in a vertical arc followed by a horizontal arc extending to a place directly over the target node and then descending vertically to that node. As the reader can see from the figure, columns in between the top and bottom nodes have horizontal arcs that can be paired with one another uniquely. Each such paring can be assigned an element $U_{k}$ and the rectilinear diagram gives a product of the $U_{k}$ that represents the element $Q.$  In fact, the form of the product of generators produced by this method is identical with the Jones normal form \cite{KD} for elements of the Temperley-Lieb Monoid. One knows that two elements in $TLM_{n}$ are algebraically equivalent if and only if they have the same normal form. This implies, by this construction that: {\it The Connection Monoid $CM_{n}$ is isomorphic as an algebraic structure with the Temperley-Lieb Monoid.} See \cite{KD,KDosen,KDosen1} for other discussions of this fact. $\hfill\Box$\\

The translation of the abstract Temperley-Lieb Monoid to the combinatorial Connection Monoid is useful for understanding many questions that we shall examine in this paper.\\

\begin{figure}[htb]
     \begin{center}
     \begin{tabular}{c}
     \includegraphics[width=7cm]{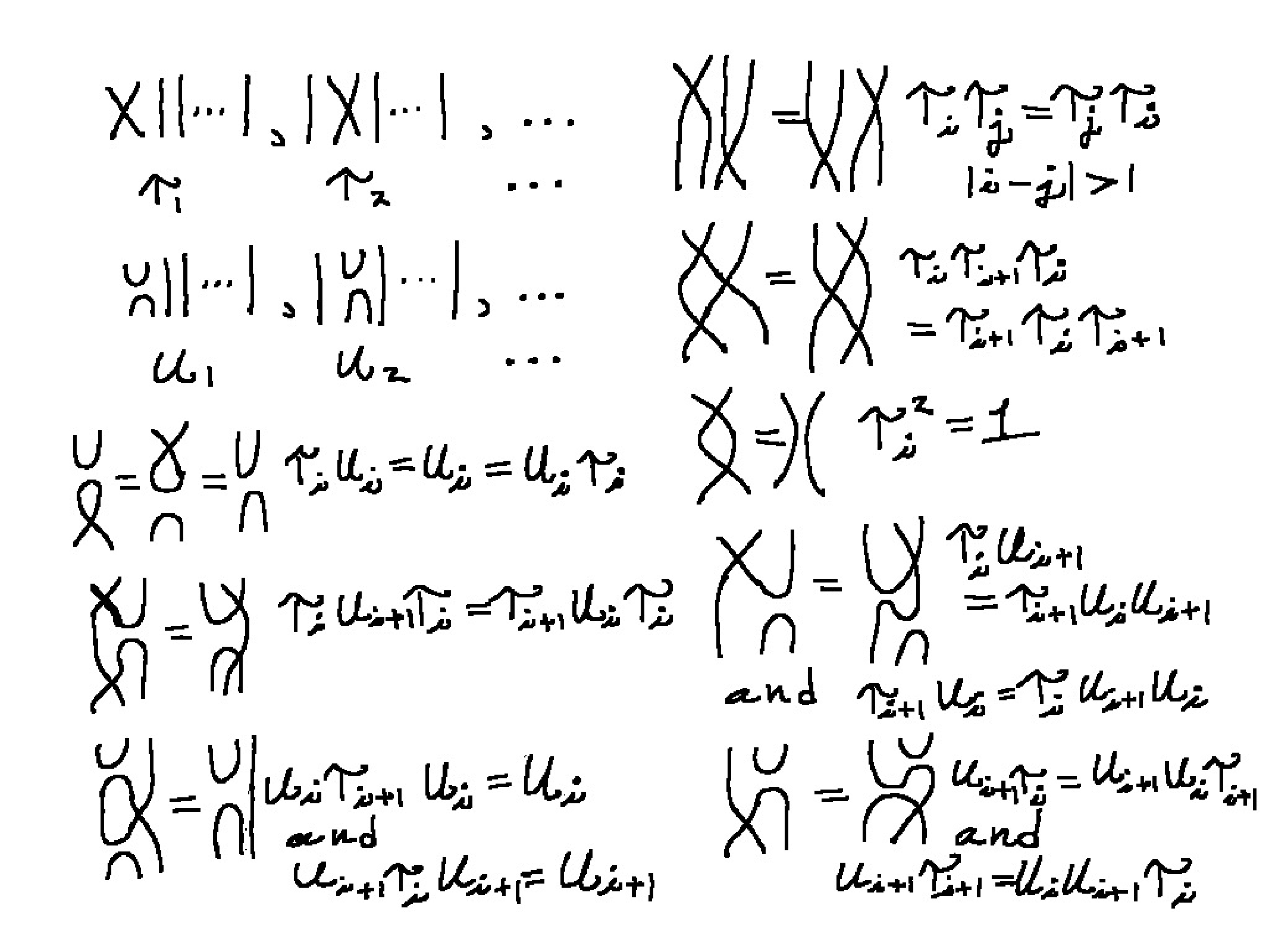}
     \end{tabular}
     \caption{\bf Brauer Monoid Relations}
     \label{Brauer2}
\end{center}
\end{figure}

\noindent {\bf Definition.} A similar description of an algebra associated with the Brauer Connection Monoid is also available. We define the {\it Brauer Monoid} on $n$ strands,  $BM_{n}$,  to be the algebraic structure with a scalar element $d$ and a unit element $1$, just as in the $TL_{n}$ and generators $U_{1},\cdots U_{n-1}$ and $\tau_{1}, \cdots \tau_{n-1}$ with the relations listed below. These relations are the Temperley-Lieb relations for the $U_i$, symmetric group relations for the $\tau_{i},$ and relations that involve the interaction of the Temperley-Lieb elements and the elements of the symmetric group. See Figure~\ref{Brauer2} for diagrammatic illustrations of these relations.
\begin{enumerate}
\item $U_i^2 = \delta U_i$ for all $i = 1, 2, \cdots , n-1.$
\item $U_i U_j U_i = U_i$ whenever $|i-j|=1.$
\item $U_i U_j = U_j U_i$ whenever $|i-j|>1.$
\item $\tau_{i}^{2} = 1$ for all $i.$
\item $\tau_i \tau_j \tau_i = \tau_j \tau_i \tau_j$ whenever $|i-j|=1.$
\item $\tau_i \tau_j = \tau_j \tau_i$ whenever $|i-j|>1.$
\item $\tau_i U_j = U_j \tau_i$ whenever $|i-j|>1.$
\item $\tau_{i} U_{i} = U_{i}\tau_{i} = U_{i}$ for all $i.$
\item $\tau_i U_{i+1} \tau_i = \tau_{i+1} U_i \tau_{i+1}$ for all $i< n-1.$
\item $U_{i} \tau_{i+1} U_{i} = U_{i}$ for all $i< n-1.$
\item $U_{i+1}\tau_{i}U_{i+1} = U_{i+1}$ for $i< n-1.$
\item $\tau_{i}U_{i+1} = \tau_{i+1}U_{i}U_{i+1}$ for all $i< n-1.$
\item $\tau_{i+1}U_{i} = \tau_{i}U_{i+1}U_{i}$ for all $i< n-1.$
\item $U_{i+1}\tau_{i} = U_{i+1}U_{i}\tau_{i+1}$ for all $i< n-1.$
\item $U_{i}\tau_{i+1} = U_{i}U_{i+1}\tau_{i}$ for all $i< n-1.$
\end{enumerate}
 
\begin{figure}[htb]
     \begin{center}
     \begin{tabular}{c}
     \includegraphics[width=7cm]{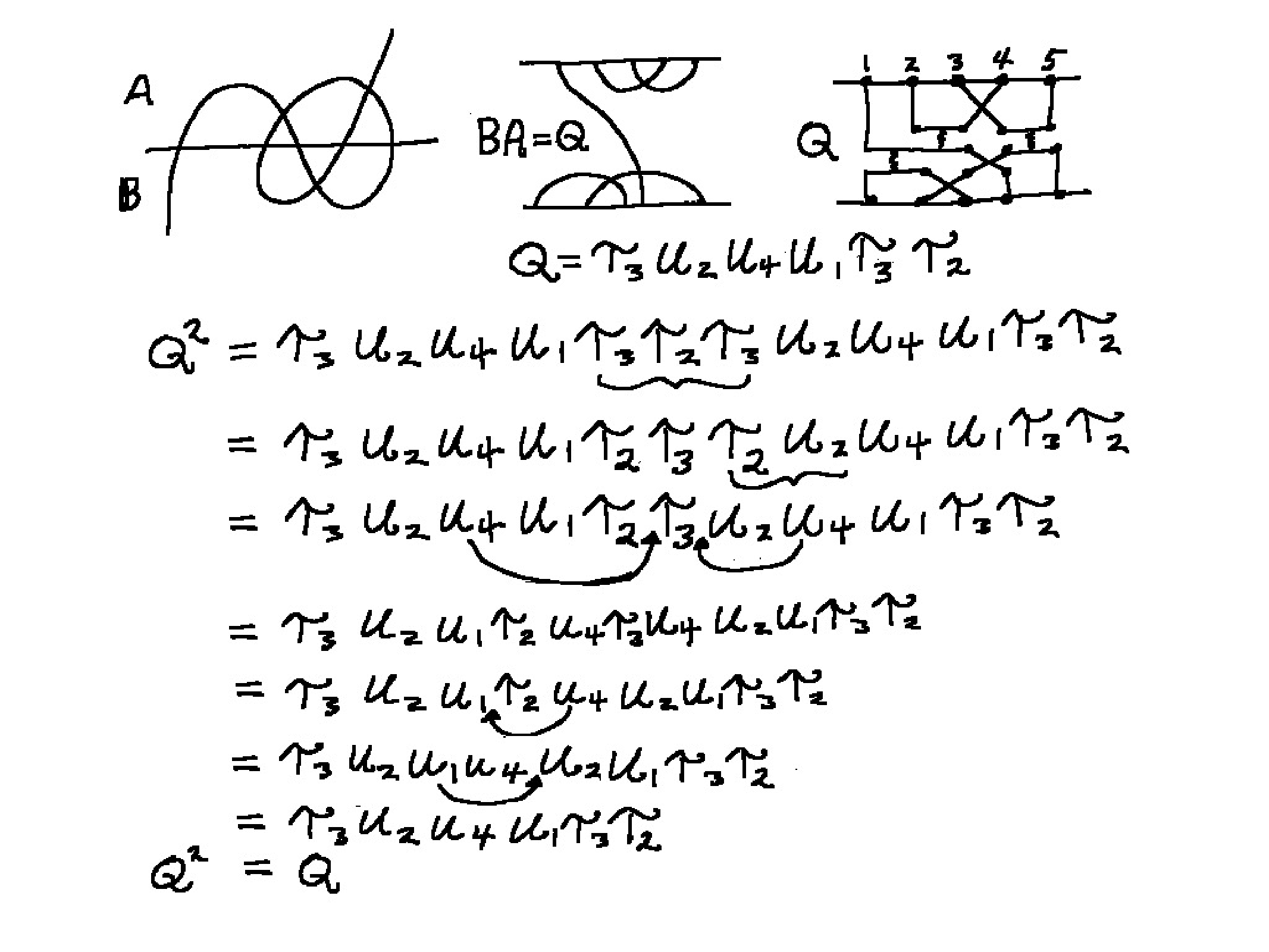}
     \end{tabular}
     \caption{\bf Rectilinear Method constructing an idempotent in the Brauer Monoid}
     \label{Brauer3}
\end{center}
\end{figure}

 Just as every diagram in  $CM_{n}$ corresponds to a canonical element in $TL_{n}$, the analog of this correspondence works for $BCM_{n}$ and $BM_{n}.$ See Figure~\ref{Brauer3} for an illustration of the correspondence. We use a generalization of the rectilinear construction described above for the Connection Monoid. In the case of the Brauer Connection Monoid we use the placement of the crossings to determine the rest of the rectilinear structure and then make pairings in the columns as before to create the $U_{i}$ factors for the algebraic expression of the connection structure. Once again it is the case that \\
 
\noindent {\bf Theorem.} {\it  Two connection structures in the Brauer Connection Monoid $BCM_{n}$  are equivalent if and only if their corresponding algebraic expressions are equivalent in the Brauer Monoid $BM_{n}.$} \\

\begin{figure}[htb]
     \begin{center}
     \begin{tabular}{c}
     \includegraphics[width=7cm]{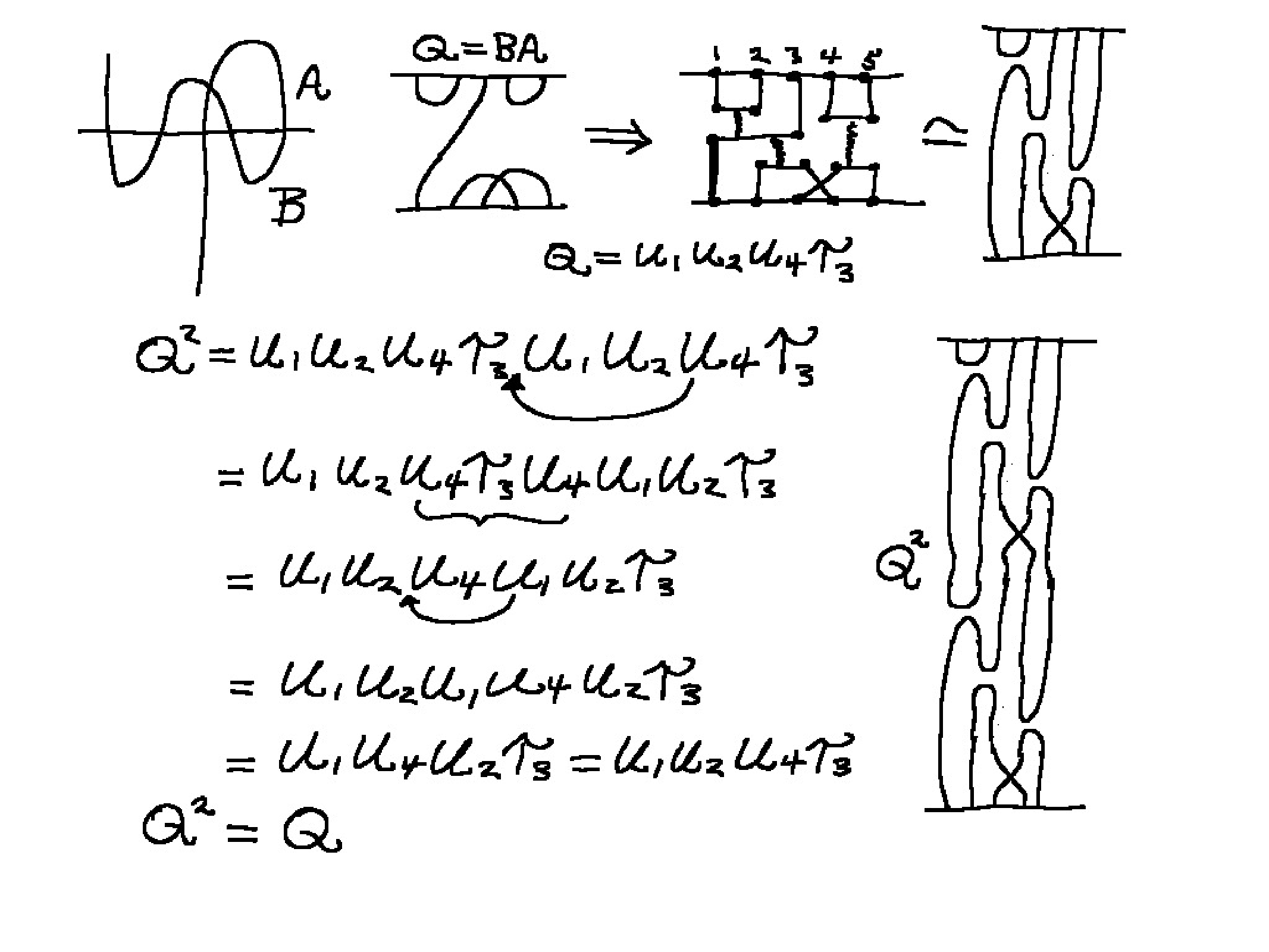}
     \end{tabular}
     \caption{\bf Rectilinear Method constructing an idempotent in the Brauer Monoid}
     \label{Brauer2.5}
\end{center}
\end{figure}

\noindent {\bf Proof.} We sketch the proof. Using \cite{Yetter} we show that two elements in the Brauer Connection Monoid are equivalent if and only if one can be obtained from the other by Morse diagrammatic moves of the type shown in 
Figure~\ref{tanglecat}, using flat crossings rather than knot theoretic crossings. The moves on algebraic elements in the Brauer Monoid are shown to produce all such categorical moves under the constraint of paired maxima and minima. The rectilinear method for converting elements of the Brauer Connection Monoid to the Brauer Monoid  (Figure~\ref{Brauer2.5})  preserves the equivalence classes and gives the desired isomorphism.  $\hfill\Box$\\
 
\begin{figure}[htb]
     \begin{center}
     \begin{tabular}{c}
     \includegraphics[width=8cm]{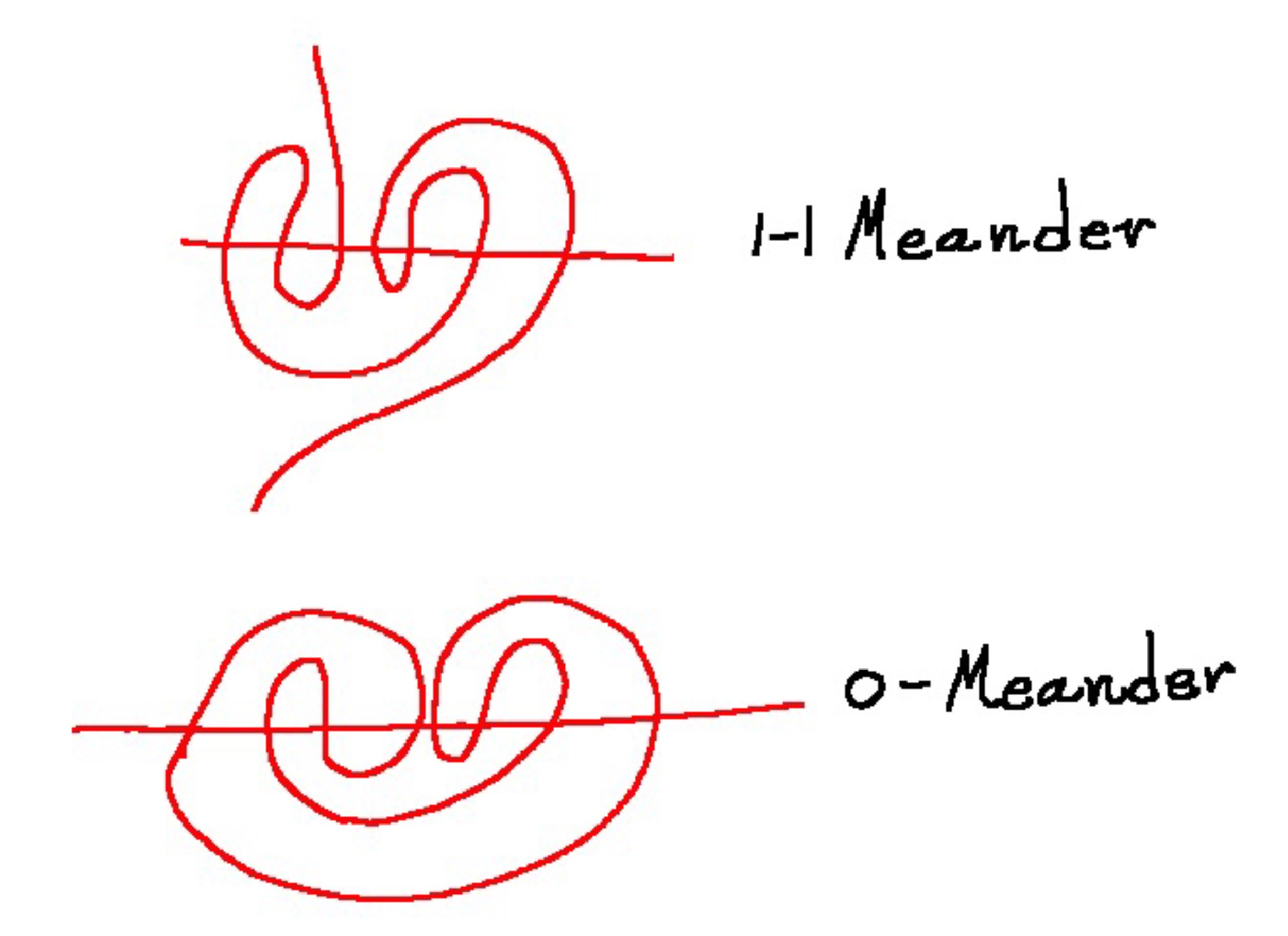}
     \end{tabular}
     \caption{\bf Meanders}
     \label{meanders}
\end{center}
\end{figure}

 \noindent {\bf Categories.} There are natural categories that generalize the Connection Monoid and the Brauer Connection Monoid. Allow two rows of nodes with different numbers of nodes in the two rows. See Figure~\ref{conncat}. Let $C_{n,m}$ denote connections with $n$ elements in the upper row and $m$ elements in the bottom row. It will be assumed that $n+m$ is even so that each node is involved in a pairing. For $C_{n,m}$ the pairings form a planar diagram. Similarly we have $BCM_{n,m}$ for morphisms in the Brauer Monoid Category. It is often the case that elements in the Connection Monoid can be factored further in the Connection Category and similarly for the Brauer Monoid and Brauer Category. Factorizations of the identity in these categories are generalizations of {\it meanders}, where a classical meander \cite{Meanders} is a Jordan curve in the plane that is partitioned by a straight line that intersects the curve transversely. See Figure~\ref{meanders}. \\
 
 When we factor the identity in $CM_{1,1}$, this is called a meander of type $(1,1)$ and can be described by starting on one side of a line and constructing a curve that goes back and forth across the line (transversely) ending up on the other side of the line. The meander can be regarded as a composition $AB$ where
 $A \in CM_{1,k}$ and $B \in CM_{k,1}$ so that $AB$ is equivalent to the identity in $CM_{1,1}.$ Then $P=BA \in CM_{k,k} = CM_{k}$ and we have $P^2 = BABA = B1A=BA=P$ so that $P$ is an idempotent in the Connection Monoid. Similar remarks apply the Brauer Connection Monoid and Category. The generalized meanders are in the $(1,1)$ case, factorizations of self-crossing curves that are equivalent to the identity under flat Reidemeister moves.\\
 
 It is of interest to see how the idempotents produced in this way can be seen as idempotents in the Temperley-Lieb Monoid and in the Brauer Monoid.  \\
 
 \noindent {\bf Definition.}  The {\it Tangle Monoid:} $Tang_{n}$ ( for $n= 1,2,3,\cdots $) and the {\it Tangle Category:} $TangCat$ are generalizations of the Brauer Monoid and Brauer Category where the crossings are replaced by knot theoretic crossings and the relations are replaced by Reidemeister Moves among the crossings and braiding moves in the monoids. We refer the reader to \cite{KD} for definitions of the Tangle Category. Here we take tangles up to all three Reidemeister moves and we include relations in the Tangle Monoid to include all the Reidemeister moves as well. The Tangle Monoid $Tang_{n}$ on $n$ strands is generated by $U_{1},U_{2},\cdots, U_{n-1}$ satisfying Temperley-Lieb relations and $\sigma_{1},\sigma_{2},\cdots, \sigma_{n-1}$ satisfying braiding relations and by relations that involve the interaction of Temperley-Lieb and braiding generators. Note that each $\sigma_{i}$ has a multiplicative inverse denoted by $\bar{\sigma_{i}}$ such that $\sigma_{i}\bar{\sigma_{i}} = \bar{\sigma_{i}}\sigma_{i} = 1$ for all $i,$ and $1$ is a multiplicative identity in the monoid. The full set of relations is listed below.
 \begin{enumerate}
\item $U_i^2 = \delta U_i$ for all $i = 1, 2, \cdots , n-1.$
\item $U_i U_j U_i = U_i$ whenever $|i-j|=1.$
\item $U_i U_j = U_j U_i$ whenever $|i-j|>1.$
\item $\sigma_i \sigma_j \sigma_i = \sigma_j \sigma_i \sigma_j$ whenever $|i-j|=1.$
\item $\sigma_i \sigma_j = \sigma_j \sigma_i$ whenever $|i-j|>1.$
\item $\sigma_i U_j = U_j \sigma_i$ whenever $|i-j|>1.$
\item $\sigma_{i} U_{i} = U_{i}\sigma_{i} = U_{i}$ for all $i.$
\item $\bar{\sigma}_{i} U_{i} = U_{i}\bar{\sigma}_{i} = U_{i}$ for all $i.$
\item $\sigma_i U_{i+1} \bar{\sigma}_{i } = \bar{\sigma}_{i+1} U_i \sigma_{i+1}$ for all $i< n-1.$
\item $\bar{\sigma}_{i} U_{i+1} \sigma_i = \sigma_{i+1} U_i  \bar{\sigma}_{i+1}$ for all $i< n-1.$
\item $U_{i} \sigma_{i+1} U_{i} = U_{i}$ for all $i< n-1.$
\item $U_{i} \bar{\sigma}_{i+1} U_{i} = U_{i}$ for all $i< n-1.$
\item $U_{i+1}\sigma_{i} U_{i+1} = U_{i+1}$ for $i< n-1.$
\item $U_{i+1}\bar{\sigma}_{i}U_{i+1} = U_{i+1}$ for $i< n-1.$
\item $\sigma_{i}U_{i+1} = \bar{\sigma}_{i+1}U_{i}U_{i+1}$ for all $i< n-1.$
\item $\sigma_{i+1}U_{i} = \bar{\sigma}_{i}U_{i+1}U_{i}$ for all $i< n-1.$
\item $U_{i+1}\sigma_{i} = U_{i+1}U_{i}\bar{\sigma}_{i+1}$ for all $i< n-1.$
\item $U_{i}\sigma_{i+1} = U_{i}U_{i+1}\bar{\sigma_{i}}$ for all $i< n-1.$
\item The roles of $\sigma$ and $\bar{\sigma}$ can be interchanged in the previous four relations.
\end{enumerate}

\begin{figure}[htb]
     \begin{center}
     \begin{tabular}{c}
     \includegraphics[width=7cm]{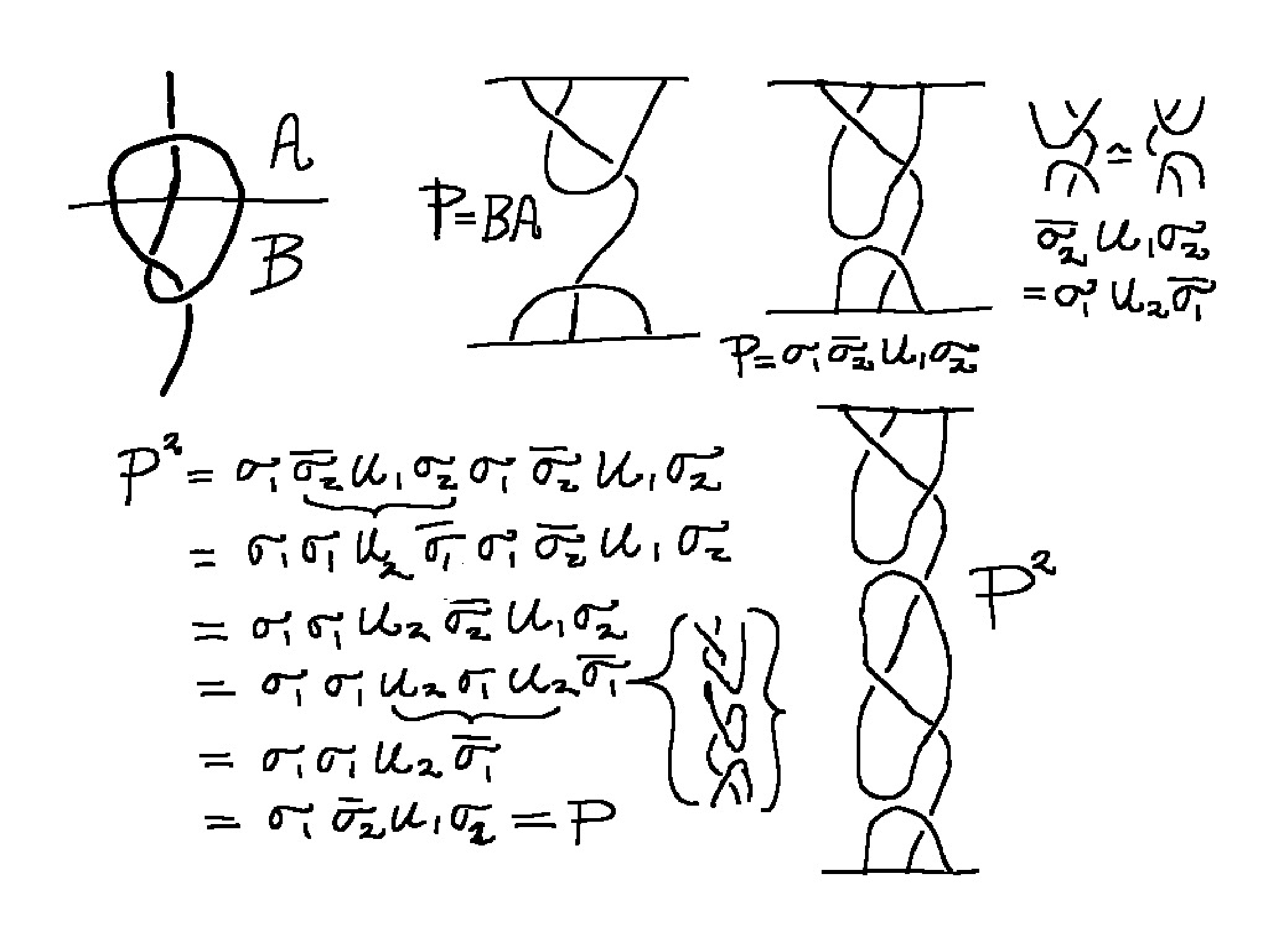}
     \end{tabular}
     \caption{\bf Constructing an idempotent in the Tangle Monoid}
     \label{tangle1}
\end{center}
\end{figure}

\begin{figure}[htb]
     \begin{center}
     \begin{tabular}{c}
     \includegraphics[width=7cm]{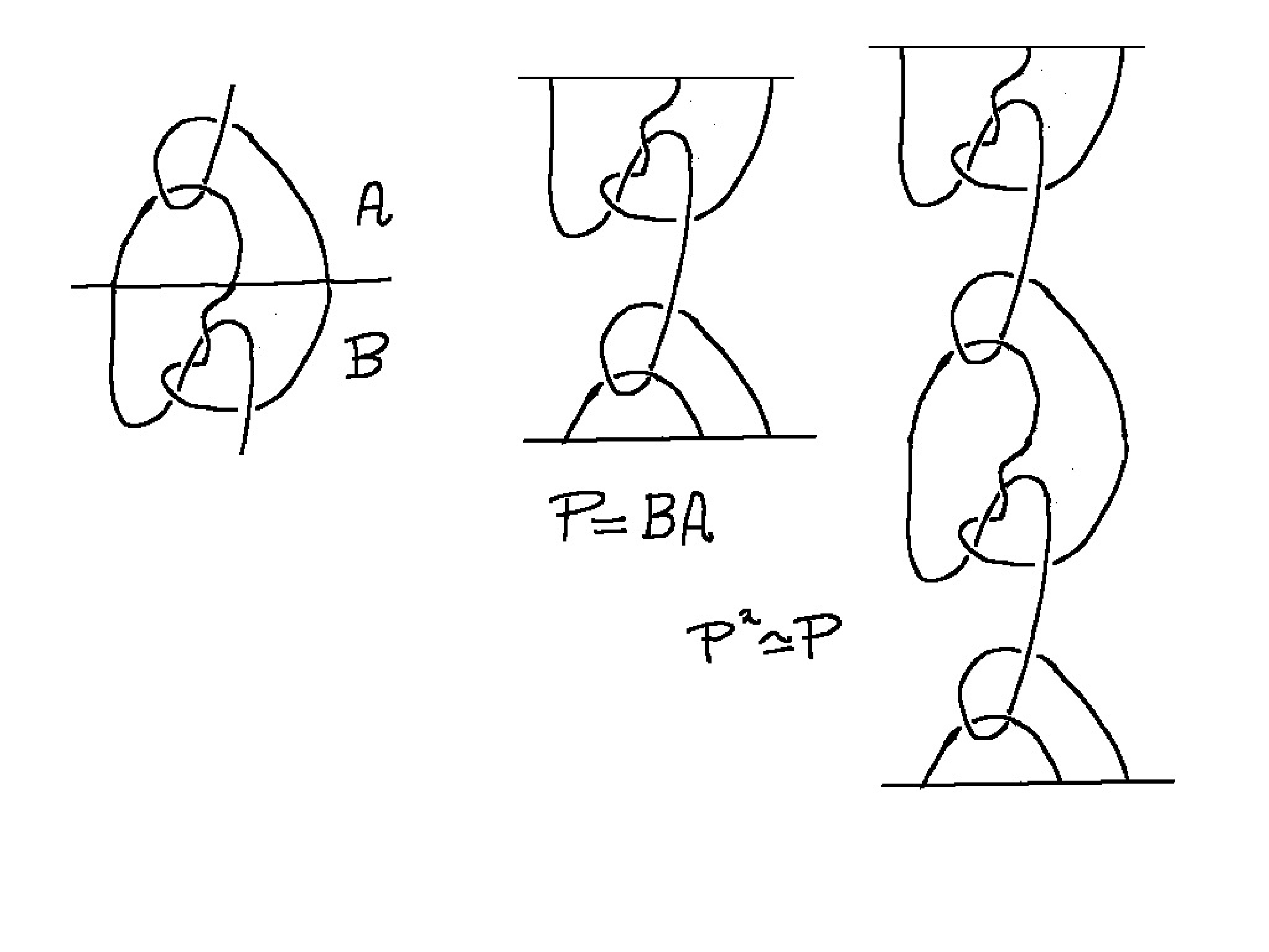}
     \end{tabular}
     \caption{\bf Using a hard unknot to construct an idempotent in the Tangle Category}
     \label{tangle2}
\end{center}
\end{figure}

\begin{figure}[htb]
     \begin{center}
     \begin{tabular}{c}
     \includegraphics[width=7cm]{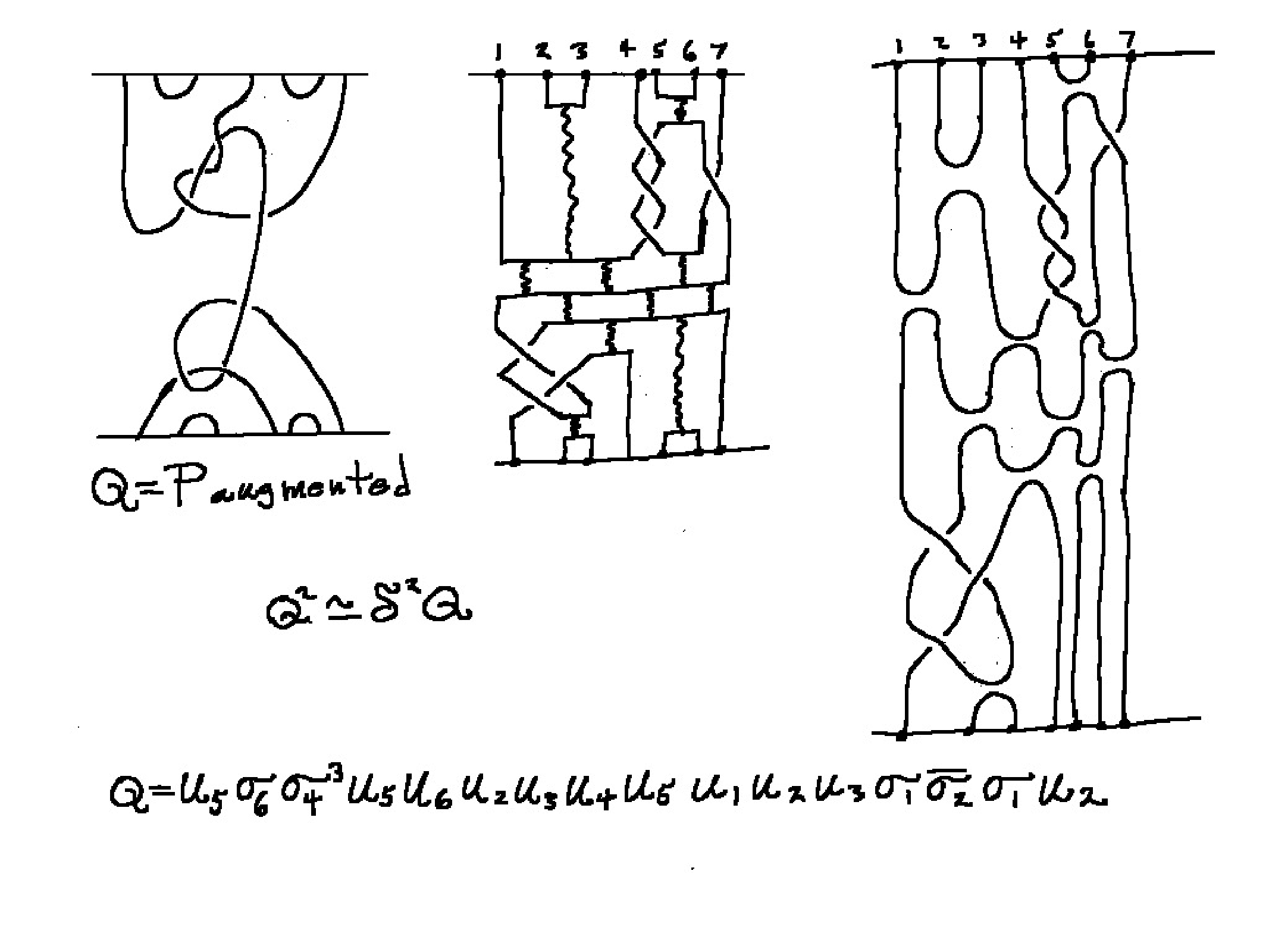}
     \end{tabular}
     \caption{\bf Embedding a hard unknot in a reduction problem in the Tangle Monoid}
     \label{tangle3}
\end{center}
\end{figure}

We can produce idempotents in the tangle category from factorizations of the identity in that category. These idempotents do not always directly yield idempotents in the Tangle Monoid due to the fact that an element in the Tangle Category does not always 
result in paired maxima and minima giving the needed Temperley-Lieb elements for a factorization into Monoid generators. See Figure~\ref{tangle1} for a successful example of this kind and see Figure~\ref{tangle2} for an example where the pairing is incomplete.
In the case of this example, we have a factorization of the identity using a ``hard unknot" \cite{HardUnknots}. This means that the unknot in question cannot be directly simplified by Reidemeister moves without first making at least one Reidemister move that increases the complexity of the diagram (complexity measured by the number of crossings in a diagram). The Figure~\ref{tangle3} shows how the rectilinear method produces a factorization for an element $Q$ in the Tangle Monoid where $Q$ is obtained by adding two new maxima and minima to the element $P$ in the Tangle Category where $AB$ is the factorization of the identity and $P=BA.$ we have in the category that $P^2 = P$ and we have in the Monoid that $Q^2 = \delta^2 Q.$ Since our construction uses a hard unknot we see that
$Q^2$ cannot be reduced directly to $\delta^2 Q$ without first performing algebraic moves that increase the complexity of its expression (where complexity is measured in terms of the number of braiding generators in the expression). This example shows how the 
topological complexity of factorizations of the identity can lead to algebraic complexity of their corresponding monoid expressions.\\

\begin{figure}[htb]
     \begin{center}
     \begin{tabular}{c}
     \includegraphics[width=8cm]{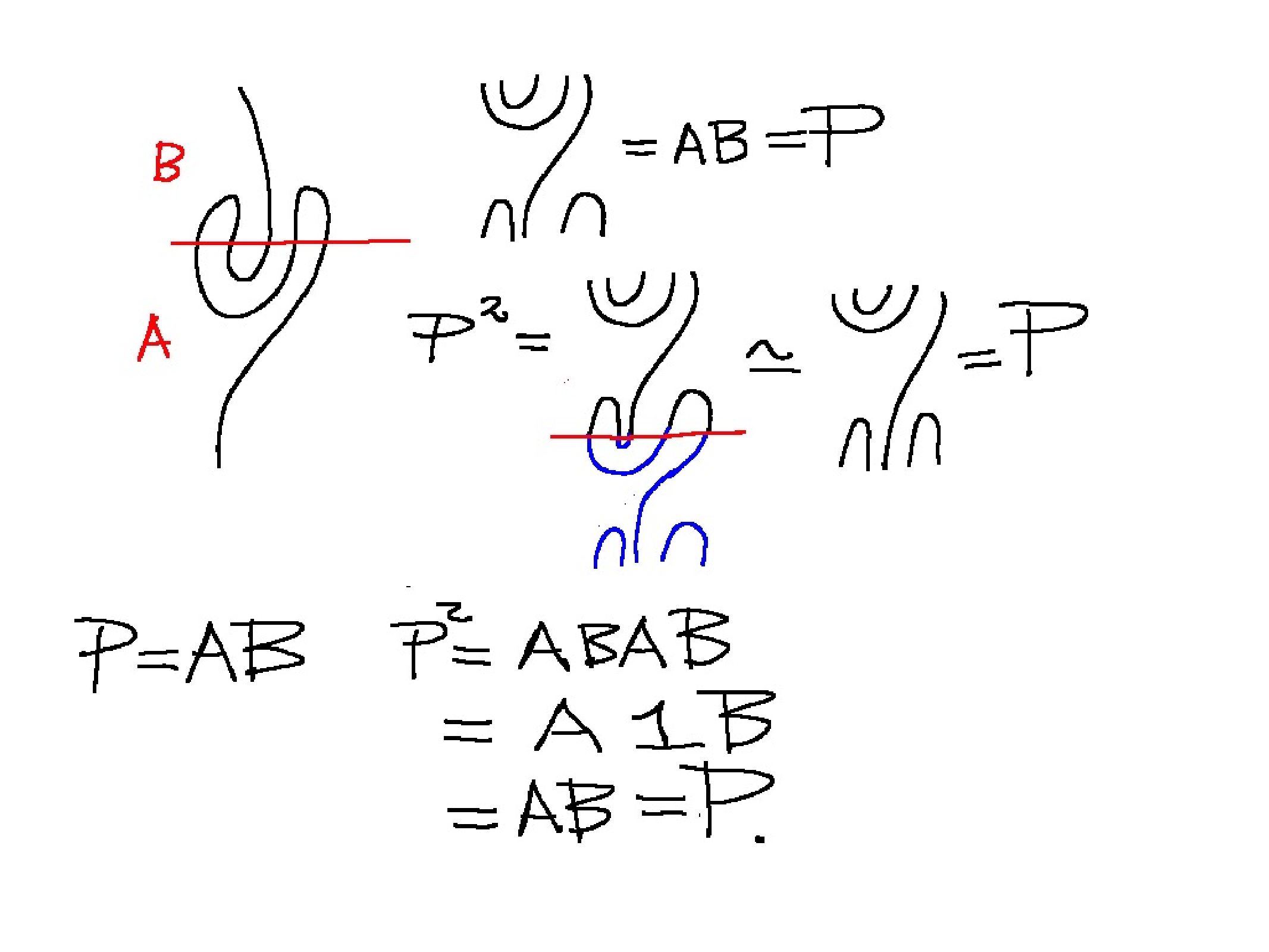}
     \end{tabular}
     \caption{\bf Factorizing an Identity}
     \label{factorid}
\end{center}
\end{figure}

\begin{figure}[htb]
     \begin{center}
     \begin{tabular}{c}
     \includegraphics[width=8cm]{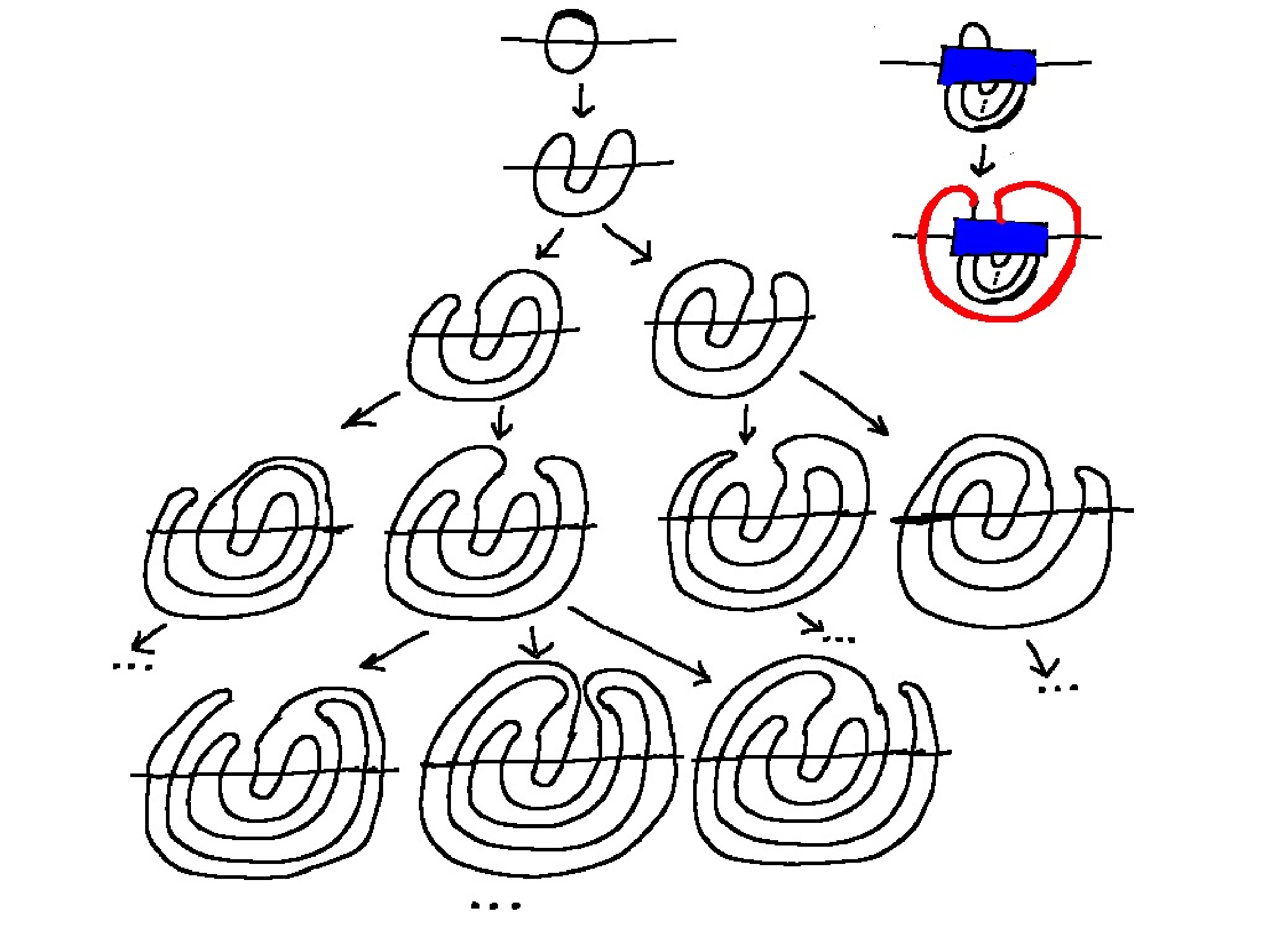}
     \end{tabular}
     \caption{\bf Enumerating Meanders}
     \label{enum}
\end{center}
\end{figure}

\begin{figure}[htb]
     \begin{center}
     \begin{tabular}{c}
     \includegraphics[width=8cm]{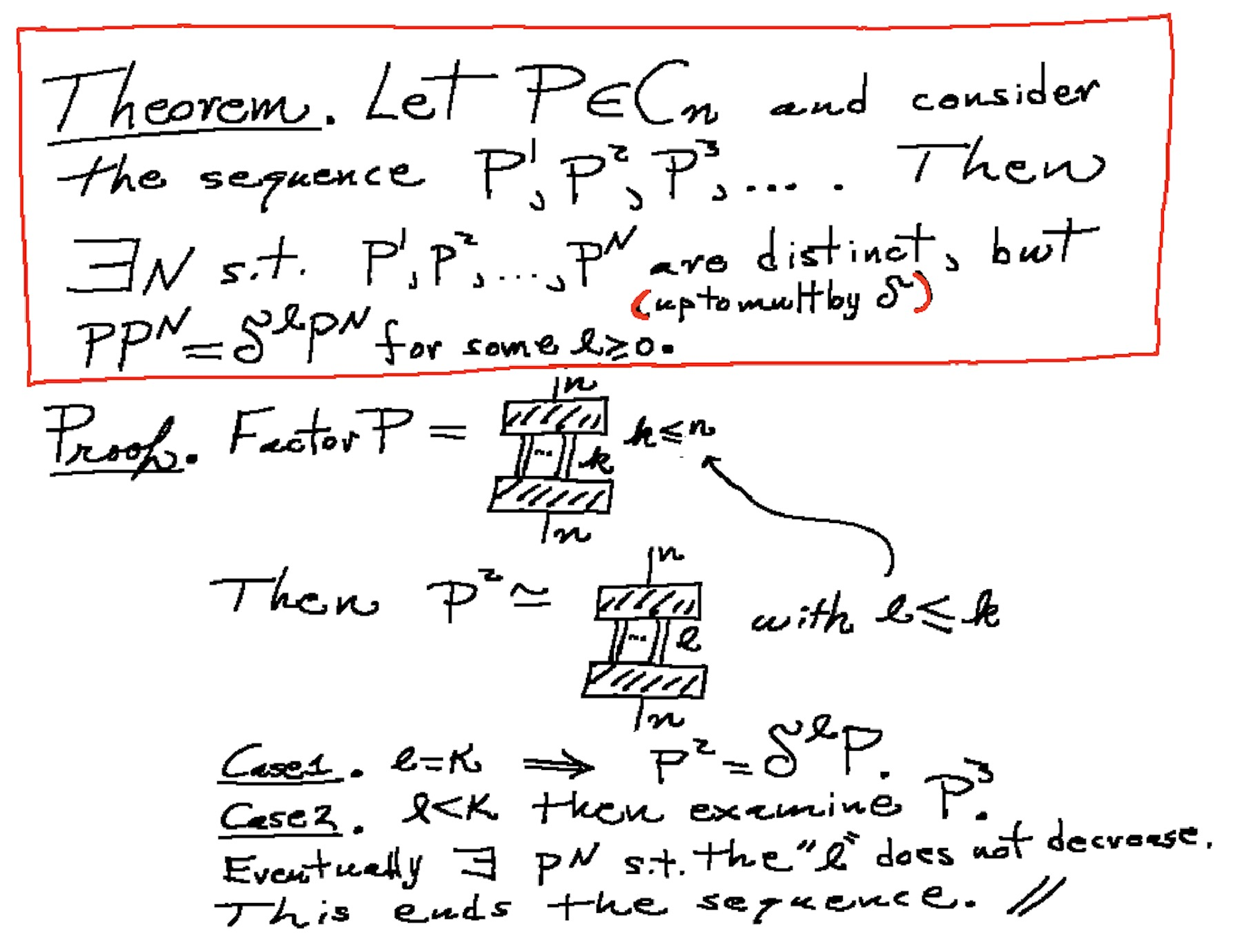}
     \end{tabular}
     \caption{\bf Powers of Elements in the Connection Monoid}
     \label{theorem}
\end{center}
\end{figure}

\begin{figure}[htb]
     \begin{center}
     \begin{tabular}{c}
     \includegraphics[width=8cm]{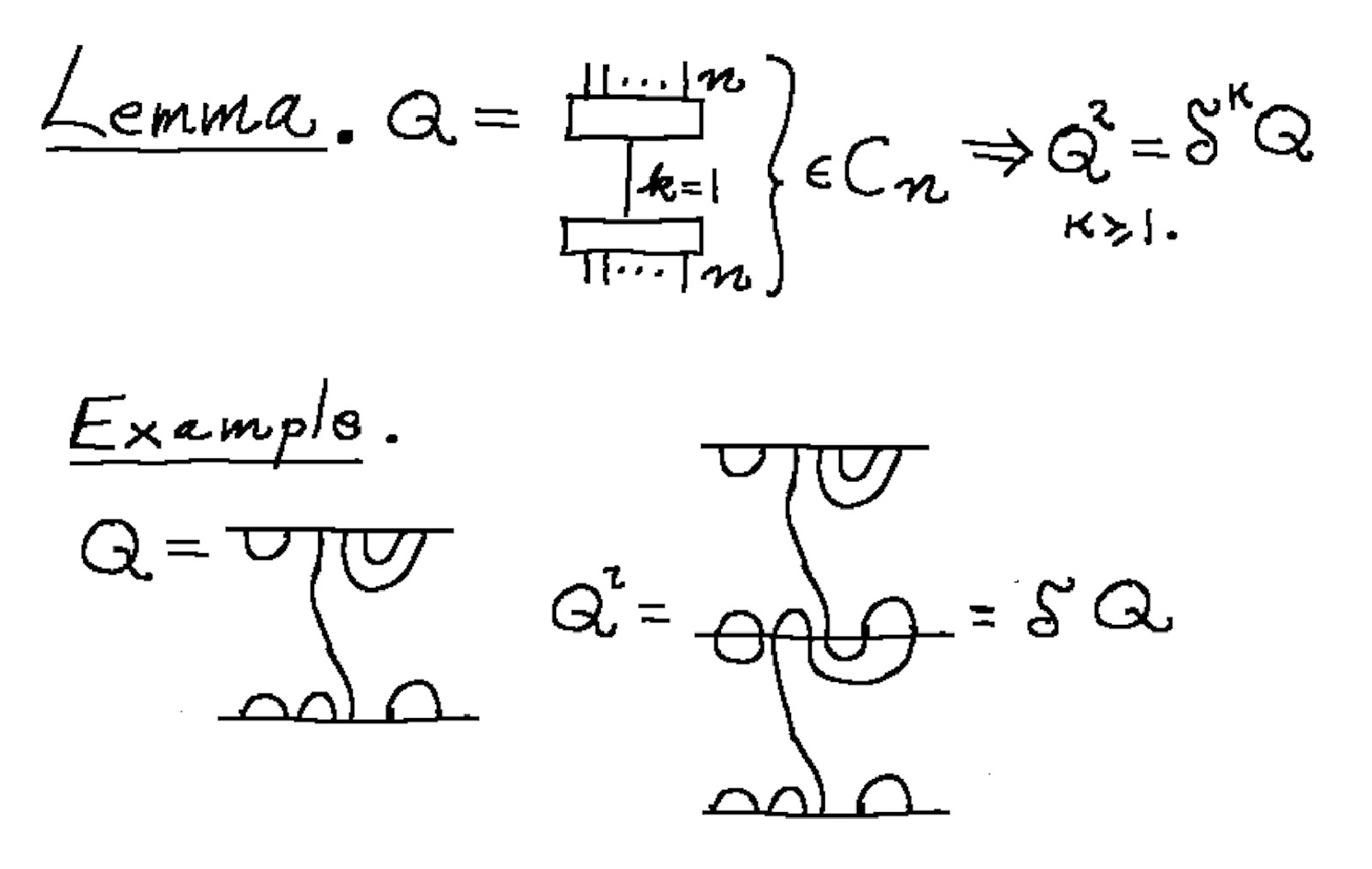}
     \end{tabular}
     \caption{\bf In the Temperley-Lieb Monoid, when k=1, then $Q^2 = \delta^{K}Q$ }
     \label{lemma}
\end{center}
\end{figure}

\begin{figure}[htb]
     \begin{center}
     \begin{tabular}{c}
     \includegraphics[width=8cm]{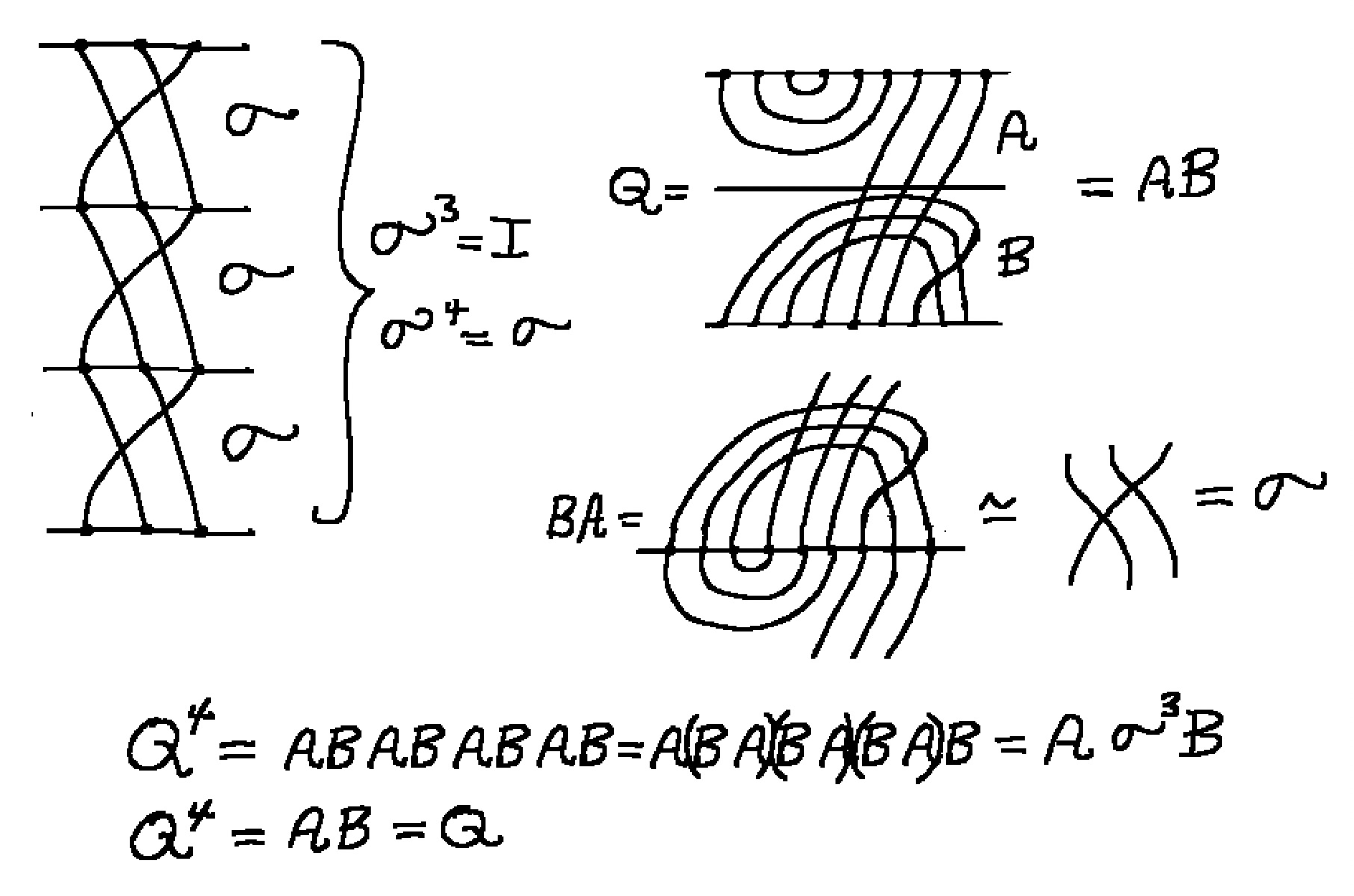}
     \end{tabular}
     \caption{\bf Brauer Examples}
     \label{brauer}
\end{center}
\end{figure}

\begin{figure}[htb]
     \begin{center}
     \begin{tabular}{c}
     \includegraphics[width=8cm]{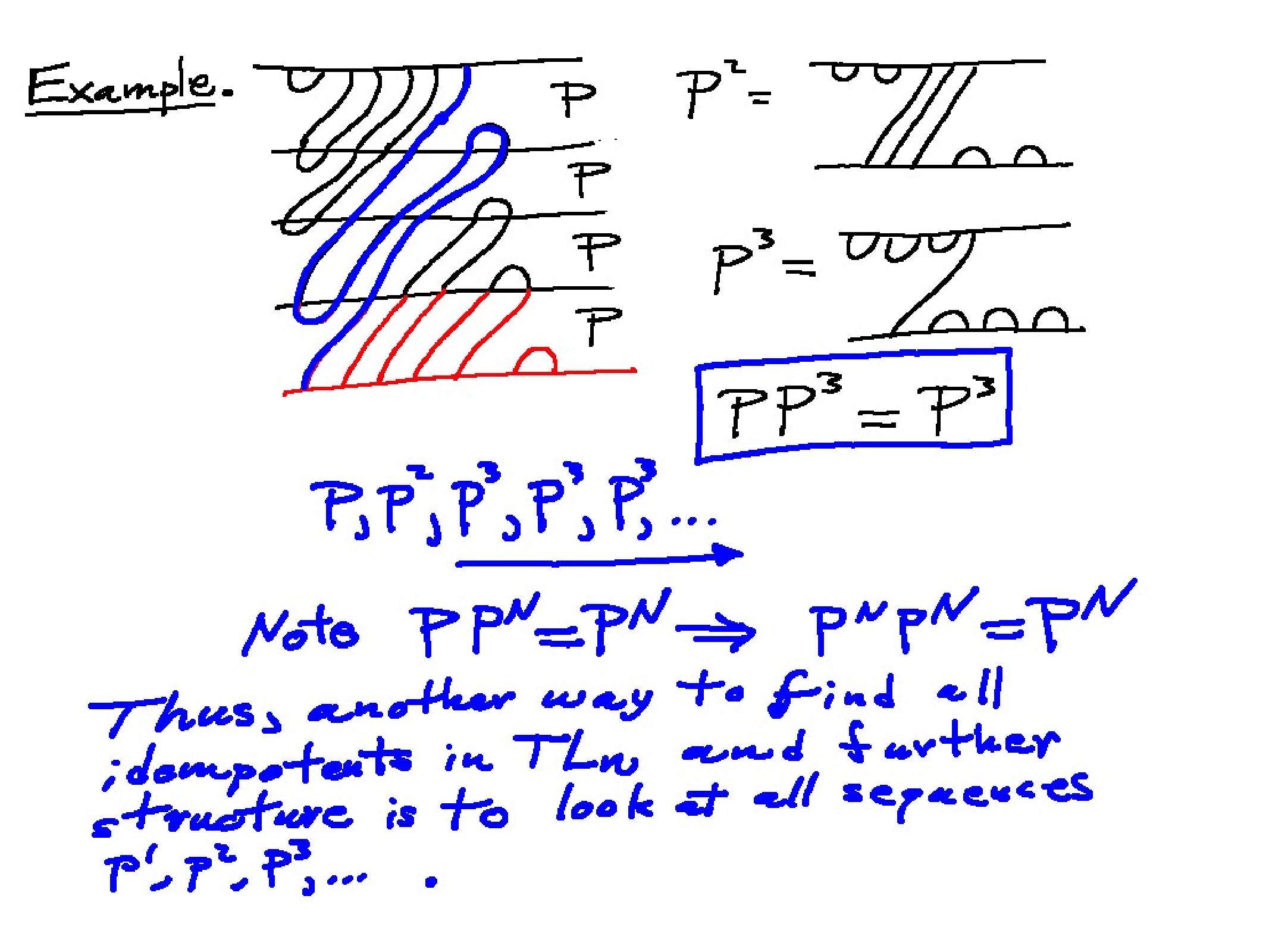}
     \end{tabular}
     \caption{\bf Powers  Example 1}
     \label{example1}
\end{center}
\end{figure}

\begin{figure}[htb]
     \begin{center}
     \begin{tabular}{c}
     \includegraphics[width=7cm]{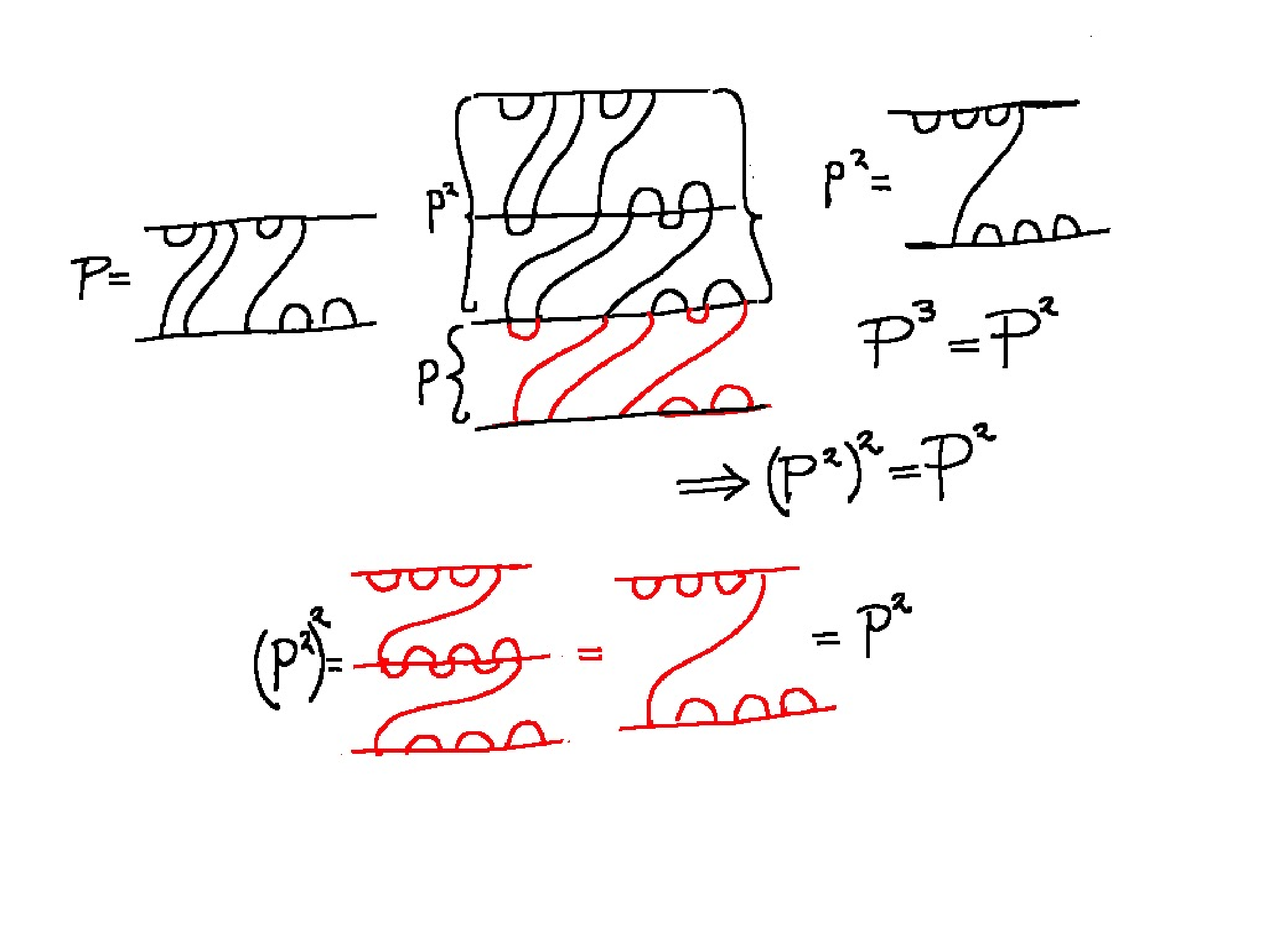}
     \end{tabular}
     \caption{\bf Powers Example 2}
     \label{example2}
\end{center}
\end{figure}

 \section{Examples and Structures}
 In this section we systematically comment on each of the Figures in the paper and what they mean in terms of the structure of idempotents in different categories and monoids.\\
 \begin{enumerate}
 \item Figure~\ref{examples}: In this figure the example $\alpha = U_1 U_2$ in the Temperley-Lieb algebra is shown. This element is idempotent since $\alpha^{2} = U_1 U_2 U_1 U_2 = U_1 U_2$ where the reduction comes from the 
 Temperley-Lieb relation $U_1 U_2 U_1 = U_1.$ The figure illustrates how diagrammatically, $\alpha = AB$ in the Temperley-Lieb Category and that $BA$ is a factorization of the identity so that $ABAB = A1B = AB.$ This figure illustrates the main ideas in the paper.
 \item Figure~\ref{TLDiagrams}:This figure illustrates for three strands, how the diagrams for the Temperley-Lieb algebra embody the relations in the algebra via topological simplification of arcs in the plane.
 \item Figure~\ref{tanglecat}: This figure illustrates the generators for the tangle category, consisting in cup, cap and crossing morphisms and an identity arc. Since this is a tensor category, all morphisms are built via composition and juxtaposition of the generating morphisms.
 \item Figure~\ref{repques}: This figure illustrates succinctly the diagrammatic version of the Jones representation of the Artin Braid Group to the Temperley-Lieb algebra and how composition of this representation with a trace function on the algebra
 (diagrammatically countling loops in the diagram closures) gives the bracket model of the Jones polynomial. A key question about this representation is whether it is faithful for all $n$ where $n$ is the number of braid strands. Understanding the multiplicative structure of the Temperley-Lieb monoid may help in understanding this question.
 \item Figure~\ref{conncat}: This figure illustrates the diagrammatics of the planar connection monoid. Every element of the Temperley-Lieb monoid has a corresponding element in the Connection Monoid and these structures are algebraically equivalent.
 \item Figure~\ref{rect1}: This figure illustrates how an element in the Connection Monoid can, by converting it to a rectilinear diagram, be mapped canonically to a specific element in the Temperley-Lieb monoid. The corresponding element in the Temperley-Lieb monoid is in Jones normal form \cite{KD,J1,J2,J3}. From this one can see that the Connection Monoid and the Temperley-Lieb Monoid are isomorphic. This figure also shows how an element in the Connection Monoid may factor into a product $BA$ in the Connection Category in such a way that $AB$ is a factorization of the identity. This means that $P = AB$ is an idempotent in the Connection Monoid. The figure illustrates how this element of the Connection Monoid can be translated (by the rectilinear method) to a specific product in the Temperley-Lieb monoid. The figure shows how one can verify directly that the resulting product is idempotent. The problem of direct characterization of idempotents in the Temperley-Lieb monoid is wide open. In general such idempotents will occur via more general factorizations of the identity (with any number of strands) in the Connection Category. To understand the idempotents at the level of the Temperley-Lieb monoid and its relations appears to be very difficult.
 \item Figure~\ref{rect2}: This figure illustrates the construction of an idempotent in the Temperley-Lieb monoid from a more complex meander (factorization of the identity) but along the lines of the previous example.
 \item Figure~\ref{factorid}: This figure illustrates a simple example producing an idempotent in the five-strand connection category.
 \item Figure~\ref{meanders}: This figure illustrates a $1-1$ meander (factorization of the identity) and a classical $0-0$ meander.
 \item Figure~\ref{enum}: This figure illustrates the enumeration of special meanders of type $0-0.$ In this case {\it special} means that the bottom half of the meander is a standard nest of cups. The method of enumeration corresponds to cutting an exterior arc from the top of a given meander, and inserting a new arc that travels across the body of the top part, crosses the horizontal line, makes a new cup at the bottom and reconnects by crossing again at the other side of the meander.
 \item Figure~\ref{theorem}: This figure illustrates the proof of the following 
 
 \noindent {\bf Theorem. } {\it Let $P \in C_n$ be an element in the Connection Monoid on $n$ strands. Consider the sequence of powers of this element $P, P^2, P^3, \cdots.$  Then there exists an integer $N$ such that $P, P^2, \cdots , P^N $ are distinct but
$P P^N = \delta^k P^N$ for some $k$ where $\delta$ is the loop value in the Monoid corresponding to the occurrence of a closed loop.}\\

\noindent {\bf Proof.} An $n$-strand element in the Connection Monoid can be factored in the Category so that there are $k$ middle strands with $k \le n.$ When an element is squared, or two elements are multiplied, then $k$ may 
decrease or remain the same. When $k$ remains the same then  $P = AB$ so that $BA$ is a $k$ to $k$ element in the category and $(BA)^{2}$ is also $k$ to $k$. The only way this can happen in the Planar Connection Monoid is if $BA = \delta^{K} I_{k}$ where $I_{k}$ denotes the identity on $k$ arcs. Thus it follows that $P^{2} = A(BA)B=\delta^{K} AI_{k}B = \delta^{K}P $ for some $K.$ This implies the Theorem since there will be some power of the initial element where the middle strand number does not decrease on further multiplication. $\hfill\Box$\\

 The next figures give examples.\\ 

\item Figure~\ref{lemma}: In this figure we illustrate the Lemma that if $P$ is in the Connection Monoid and $P$ factors in the Category so that there is a single connecting strand, then $P^2 = \delta^{K}P$ for some non-negative integer $K.$ It is interesting
to realize this fact in relation to our previous remarks about factorizing the identity in a category.\\
\item Figure~\ref{brauer}: In fact, the Theorem of Figure~\ref{theorem} is true if we replace Connection Monoid by Brauer Connection Monoid but the specific phenomena are different. Since crossed arcs can produce permutations, one can easily write elements $\sigma$ such that $\sigma^{N} = I$ and $\sigma$ has order $N.$ Then $\sigma^{N+1} = \sigma$ satisfies the pattern of the Theorem. In Figure~\ref{brauer} we illustrate an example $Q$ where a permutation is coded in one of the factors in the Category, but appears as a permutation factor when we form $Q^2.$ In particular, $Q=AB$ where $BA = \sigma$ and $\sigma$ is a permutation with $\sigma^3 = I$ where $I$ denotes the identity permutation. Thus
$$Q^4 = ABABABAB = A(BA)(BA)(BA)B = A \sigma^3 B = AB = Q.$$ The difference in the proof of the Theorem above, is that when an element $Q$ factors as $Q=AB$ with $BA$ on $k$ strands to $k$ strands, then $BA$ can contain a permutation of some finite order and still have the property that all its powers are on $k$ strands. Nevertheless, the permutation is of some finite order and so there will be an $N$ so that $(BA)^{N} = \delta^{K} I$ for some positive integer $K.$ Then we have
$Q^{N+1} = A(BA)^{N} B = \delta^{K} Q.$ This is the main case in proving the \\

\noindent {Brauer Monoid Theorem.} {\it Let $Q \in BCM_n$ be an element in the Brauer Connection Monoid on $n$ strands. Consider the sequence of powers of this element $Q, Q^2, Q^3, \cdots.$  Then there exists an integer $N$ such that $Q, Q^2, \cdots , Q^N $ are distinct but $Q Q^N = \delta^{K} Q^N$ for some $K$ where $\delta$ is the loop value in the Monoid corresponding to the occurrence of a closed loop.}\\

\noindent {\bf Proof.} The proof is indicated prior to the statement of the Theorem. $\hfill\Box$\\

 \item Figure~\ref{example1}:This figure illustrates and element $P$ in the Connection Monoid such that $P,P^2,P^3$ are distinct but $P P^3 = P^3.$ 
 \item Figure~\ref{example2}: This figure illustrates and element $P$ so that $P$ and $P^2$ are distinct but $P P^2 = P^2.$ Thus $P^2 P^2 = P P P^2 = P P^2 = P^2$ so that $P^2$ is an idempotent, and makes it apparent that it comes from a simple meander.
 \item Figure~\ref{example3}: This figure illustrates the Theorem in the case of a family of elements that Morrison \cite{Morrison} uses to construct Jones-Wenzl projectors. The Theorem may shed light on properties of these projectors.
 \item Figure~\ref{Brauer1}: This figure shows a factorization of the identity and the construction of a corresponding idempotent in the Brauer Monoid. The Brauer Monoid generalizes the Connection Monoid by allowing all possible pairings and hence crossing arcs in the planar representations.
 \item Figure~\ref{Brauer2}: This figure illustrates the basic relations in the Brauer Monoid generated by crossovers $\tau_i$ and Temperley-Lieb generators $U_i.$
 \item Figure~\ref{Brauer2.5}: This figure illustrates a factorization of the identity (generalized meander) in the Brauer Category and the construction and explicit factorization for this idempotent in the Brauer Monoid. The rectilinear construction is
 used in this context once again.
 \item Figure~\ref{Brauer3}: This figure illustrates another generalized meander and its corresponding idempotent in the Brauer Monoid.
 \item Figure~\ref{tangle1}: In this figure we factor the identity as a $1-1$ element in the Tangle Category. In this case the corresponding idempotent is directly represented in the Tangle Monoid.
 \item Figure~\ref{tangle2}: In this figure we illustrate a factorization of the identity in the Tangle Category that is itself a ``hard unknot". There is no sequence of simplifying Reidemeister moves that undoes this knot. This means that in the category the representing diagram will have to made more complex (more crossings) before it can be undone. We show the corresponding idempotent in the Tangle Category.
 \item Figure~\ref{tangle3}: In this figure we illustrate that the idempotent constructed in the previous example cannot be converted directly to an element in the Tangle Monoid due to lack of pairings of cups and caps. We have augmented the 
 idempotent $P$ to a new element $Q$ with two new cups and two new caps and we show, by the rectangular technique, how to write $Q$ in the Tangle Monoid. Because $Q$ is derived from a hard unknot, $Q^2$ will have to increase its number of crossing generators in order to eventually simplify to $\delta^2 Q.$ This example shows how complexities of knotting and unknotting are reflected in properties of the algebra of the Tangle Monoid.
 \end{enumerate}

\end{document}